\documentclass[journal,subeqn,caption]{IEEEtran}

\usepackage{mathtools}
\usepackage{cite}
\usepackage{float}
\usepackage{multicol}
\usepackage{multirow}
\usepackage{hyperref}
\usepackage{graphicx}
\usepackage{amsfonts}
\usepackage{amsmath}
\usepackage{color}
\usepackage{epstopdf}
\usepackage{caption}
\usepackage{subcaption}
\usepackage{mathtools, nccmath}
\usepackage{multirow}
\usepackage{bbm}
\usepackage{marginnote}
\usepackage{tikz}
\usepackage{pbox}
\usepackage{xcolor}
\usepackage{comment}
\usepackage{amsthm}
\usepackage{url}

\usepackage{cleveref}
\crefformat{equation}{\textup{#2(#1)#3}}
\crefrangeformat{equation}{\textup{#3(#1)#4--#5(#2)#6}}
\crefmultiformat{equation}{\textup{#2(#1)#3}}{ and \textup{#2(#1)#3}}
{, \textup{#2(#1)#3}}{, and \textup{#2(#1)#3}}
\crefrangemultiformat{equation}{\textup{#3(#1)#4--#5(#2)#6}}%
{ and \textup{#3(#1)#4--#5(#2)#6}}{, \textup{#3(#1)#4--#5(#2)#6}}{, and \textup{#3(#1)#4--#5(#2)#6}}

\usepackage{algorithm}
\usepackage{algorithmic}

\DeclareMathOperator{\st}{s.t.} 
\DeclareMathOperator{\lag}{\mathcal{L}}
\DeclareMathOperator{\B}{\mathcal{B}}
\DeclareMathOperator{\G}{\mathcal{G}}
\DeclareMathOperator{\LB}{LB}
\DeclareMathOperator{\UB}{UB}
\DeclareMathOperator{\proba}{\mathbb{P}}
\DeclareMathOperator{\E}{\mathbb{E}}
\DeclareMathOperator{\cost}{Cost}
\DeclareMathOperator{\V}{\mathbb{V}}

\newcommand{\blue}[1]{\textcolor{black}{#1}}

\newcommand \dkm[1] {{\color{black}#1}}

\begin{document}
%
\title{DC Optimal Power Flow with\\ Joint Chance Constraints}
%
%

\author{Alejandra~Pe\~na-Ordieres,
        Daniel~K.~Molzahn,
        Line~A.~Roald,
        and~Andreas~W\"achter
\thanks{A. Pe\~na-Ordieres and A. W\"achter are with the Department of Industrial Engineering and Management Sciences, Northwestern University, Evanston, IL, 60208.}
\thanks{D. Molzahn is with the School of Electrical and Computer Engineering, Georgia Institute of Technology, Atlanta, GA, 30332.}
\thanks{L. Roald is with the Department of Electrical and Computer Engineering, University of Wisconsin--Madison, Madison, WI, 53706.}
}


\maketitle

\begin{abstract}
Managing uncertainty and variability in power injections has become a major concern for power system operators due to increasing levels of fluctuating renewable energy connected to the grid. This work addresses this uncertainty via a joint chance-constrained formulation of the DC optimal power flow (OPF) problem, which satisfies \emph{all} the constraints \emph{jointly} with a pre-determined probability. The few existing approaches for solving joint chance-constrained OPF problems are typically either computationally intractable for large-scale problems or give overly conservative solutions that satisfy the constraints far more often than required, resulting in excessively costly operation. 
This paper proposes an algorithm for solving joint chance-constrained DC OPF problems by adopting an S$\ell_1$QP-type trust-region algorithm. This algorithm uses a sample-based approach that avoids making strong assumptions on the distribution of the uncertainties, scales favorably to large problems, and can be tuned to obtain less conservative results. We illustrate the performance of our method using several IEEE test cases. The results demonstrate the proposed algorithm's advantages in computational times and limited conservativeness of the solutions relative to other joint chance-constrained DC OPF algorithms.

\end{abstract}

\begin{IEEEkeywords}
joint chance constraints, nonlinear optimization, optimal power flow, sample average approximation
\end{IEEEkeywords}

\IEEEpeerreviewmaketitle

\section{Introduction}

\IEEEPARstart{O}{ptimal} power flow (OPF) is a fundamental problem in power systems operations that is used for real-time operations, markets, long-term planning, and many other applications. In its classical form, OPF determines the minimum cost generation dispatch that satisfies the demand for power while adhering to network constraints and engineering limits. 

Growing quantities of renewable energy are increasing the variability and uncertainty inherent to power system operations. Many new methods account for and mitigate this uncertainty and variability~\cite{PancEtAl14}, including two- and multi-stage stochastic programming \cite{bouffard2008, morales2009, papavasiliou2013}, robust and worst-case optimization \cite{panciatici2010, jabr2013, warrington2013, lorca2015}, and chance constraints
\cite{zhang2011,VrakEtAl13,RoalEtAl13,BienCherHarn14,summers2014, BakeBern19}. 
These methods attempt to ensure secure and economical operations despite power injection uncertainty. Defining ``security'' is an important modelling question that dictates the formulation and solution algorithm. For example, robust optimization defines ``secure'' as ensuring feasibility for all realizations within a pre-specified uncertainty set, while chance-constrained optimization seeks to satisfy the constraints with a high probability $1-\alpha$, where $\alpha$ is a specified acceptable violation probability. 

We propose a formulation and solution algorithm to solve OPF problems with \emph{joint chance constraints} (JCC), which require that \emph{all} engineering limits, including both generation and line flow constraints, are satisfied simultaneously with probability $1-\alpha$. This contrasts with formulations based on single chance constraints (SCC), which split the line flow and generation limits into separate chance constraints (for each line and generator) with individual risk levels, $1-\alpha_j$, for each of those constraints. 
Allocation of risk to individual components is more straightforward in problems with SCCs, while JCCs give much stronger guarantees on overall system security. 
Generally, SCCs are much easier to solve \cite{GengXie19_I}. For example, linear SCCs with elliptical symmetric uncertainty distributions can be expressed as second-order cone programs that can be efficiently solved \cite[Lemma 2.2]{Henr07}.

Most chance-constrained OPF formulations have considered SCCs (e.g., \cite{zhang2011,RoalEtAl13,BienCherHarn14}), while a limited number have attempted to solve JCC formulations \cite{VrakEtAl13}. In \cite{BakeBern19}, a JCC problem is solved by decomposing the JCC into SCCs, which is challenging due to the difficulty in selecting the risk level for each individual constraint. Usually, the Boole or Bonferroni inequality is used to approximate the JCC. Reference \cite{ChenEtAl10} observes that even if the individual risk levels are selected optimally, the solution obtained from the SCC formulation can be suboptimal. Some efforts have been made to reduce the conservativeness of using Boole's inequality (e.g., \cite{BakeBern19}), and it has been observed that the SCC formulation leads to a low joint violation probability due to the structure of the OPF problem~\cite{RoalEtAl16}. However, in general, the SCC formulation has the following drawbacks: (1) enforcing the chance constraints individually does not give strong guarantees on the feasibility probability of the entire system, and (2) solutions that are adapted to guarantee joint feasibility can be overly conservative and costly. 

The most common methods for directly solving JCCs are based on scenario approximation (SA) (e.g., \cite{CalaCamp05,CalaCamp06,NemiShap06_SA}), which has been applied to the OPF problem in, e.g., \cite{VrakEtAl13}, and mixed-integer programming (MIP) (e.g., \cite{LuedAhme08}). Both the SA and the MIP methods provide guarantees on the quality of the solution and are sample-based approximations, meaning that they do not make assumptions on the uncertainty distributions. However, solutions from SA are often highly conservative with much lower violation probabilities than what would be acceptable and, consequently, these solutions are more costly~\cite{VrakEtAl13, roald2017chance, modarresi2018scenario}. While MIP methods converge to the desired solution with increasing sample size, the complexity of the algorithm also increases, which can result in intractability.

\dkm{This paper's main contribution is a joint chance-constrained formulation and algorithm to solve the DC~OPF problem. The formulation is based on a sample average approximation (SAA) which gives rise to a continuous non-linear programming (NLP) problem. The algorithm is an adaptation of the JCC algorithm presented in~\cite{PenaLuedWaec19}. Making this algorithm applicable to electric power systems requires careful consideration of the formulation and several non-trivial modifications to the algorithm. To improve computational performance, we (1)~select a DC~OPF formulation that avoids the need to replicate certain variables for each scenario, (2)~approximate the Hessian of the quantile function defined in~\cite{PenaLuedWaec19} in order to solve a more tractable convex problem at each iteration of the algorithm, and (3)~develop a lazy constraint generation algorithm to exploit the observation that a small number of line flow limits are binding at the solutions to typical DC~OPF problems~\cite{molzahn_roald-redundant2019}. To improve the solution quality relative to a na\"ieve application of~\cite{PenaLuedWaec19}, we propose methods for adaptively tuning the key parameters introduced in the quantile approximation. This tuning improves the out-of-sample feasibility of the resulting solutions.}



The NLP approach has several advantages:
\emph{(1)~Sample-based:} Similar to the SA and MIP methods, we use a sample-based approach that does not rely on distributional assumptions.
\emph{(2)~Scalable:} The method is scalable to large systems with many uncertain power injections where SA may be impractical due to the need for a very large sample size and MIP methods may be numerically intractable due to  the introduction of binary variables.
\emph{(3)~Tunable:} \blue{The chance constraint approximation presented in \cite{PenaLuedWaec19} depends on certain parameters that impact the conservativeness of the solution. We propose two methods that adaptively and automatically tune these parameters such that the resulting solution accurately satisfies the prescribed probability. Hence, the proposed method does not render an excessively conservative feasible region, which is an advantage over the SA method.}

The remainder of this paper is organized as follows. Section~\ref{sec: JCCP-OPF} describes the JCC-OPF formulation. Section~\ref{sec: Quantile as CC} proposes a smooth sample-based approximation of the probabilistic constraint in the JCC-OPF formulation. Section~\ref{sec: SolAlg} presents our solution algorithm. Section~\ref{sec: SelectingParameters} discusses the tuning parameters. Section~\ref{sec: NumRes} numerically demonstrates our method, benchmarked against SA. Section~\ref{sec: Conclusions} concludes the paper.

\section{Joint Chance-Constrained Optimal Power Flow} \label{sec: JCCP-OPF}
We aim to minimize the expected generation cost while satisfying all engineering limits with a high probability via a \emph{joint chance-constrained} OPF problem (JCC-OPF).
The user expresses an acceptable risk as the joint violation probability, i.e., the probability that any of the constraints are violated. 
This section formulates the JCC-OPF. This formulation is closely related to those previously presented in \cite{VrakEtAl13, RoalEtAl13, BienCherHarn14}, but differs in the handling of the forecasted operating point.

\subsubsection{Notation} Consider a power system where the sets of buses, lines, and generators are denoted by $\mathcal{B}$, $\mathcal{L}$, and $\mathcal{G}$, respectively. To simplify notation, we assume that there is one generator with active power generation $g(\omega)$ and one uncertain load $d(\omega)$, where $\omega$ is a random variable, at every bus. Then, $g(\omega), d(\omega) \in \mathbb{R}^{|\mathcal{B}|}$. If a bus $i$ does not have a generator or load, we set $g_i(\omega) = 0$ or $d_i(\omega) = 0$, respectively, whereas multiple loads or generators are handled through summation.

We use the linearized DC approximation of the active power flows \dkm{which makes the following assumptions that are standard to all DC power flow formulations~\cite{stott2009}}: 
(1)~all voltage magnitudes are $1$~per unit, (2)~neighboring buses have small angle differences, and (3)~the system is lossless. 

\subsubsection{Uncertain loads} All uncertain loads can be represented as $d(\omega) = d + \omega$, where $\omega$ is a random variable with zero mean; this can be interpreted as the sum of the forecasted value $d$ and its fluctuation $\omega$.  Due to the nature of the renewable energy uncertainty, we model $\omega$ as a continuous random variable. \dkm{We note that the random variable $\omega$ models the uncertainty in the net load, i.e., the load demands minus the outputs of stochastic renewable generators.}

\subsubsection{Generators} We model the active power generation $g(\omega)$ using an affine control policy, resembling the actions of the automatic generation control (AGC) \cite{VrakEtAl13}. Each generator adjusts its output to satisfy a fraction of the total load imbalance,
\begin{align}
    g_i(\omega) = g_i - \beta_i \Omega, \quad \forall i \in \G, \label{eq: AffinePolicy}
\end{align}
where $\Omega = \sum_{i \in \B } \omega_i$ and $\beta_i$ is the so-called participation factor of generator $i$. Our formulation's optimization variables include the generation $g$ and the participation factors~$\beta$. 

\subsubsection{Power Balance} With the lossless system representation, maintaining power balance is equivalent to ensuring that the total power generation equals the total demand,
\begin{align}
    \sum_{i\in\mathcal{G}}g_i(\omega)+\sum_{j\in\mathcal{B}}d_j(\omega)=0, \quad \forall \omega. \label{eq:balance}
\end{align}
By substituting the expressions for $g(\omega)$ and $d(\omega)$ from above, we observe that \eqref{eq:balance} is equivalent to enforcing
$$\sum_{i\in\mathcal{G}}g_i+\sum_{j\in\mathcal{B}}d_j=0 \quad\text{and}\quad \sum_{i\in\mathcal{G}} \beta_i=1.$$
Here, the first equation guarantees power balance without fluctuations $\omega=0$, while the second equation ensures system balance during fluctuations $\omega\neq 0$.

\blue{If $g(\omega)$ did not follow the affine functional form in \eqref{eq: AffinePolicy}, a copy of the variable $g(\omega)$ would need to be introduced for each realization of $\omega$ to satisfy~\eqref{eq:balance}. For a sample-based approach, this implies that the number of variables in the problem would depend on the number of scenarios considered in the sample. Thus, the choice of the functional form of $g(\omega)$ is relevant to decreasing the complexity of the JCC-OPF problem.}

\subsubsection{Power flows} We denote the line connecting buses $i$ and $j$ as $ij \in \mathcal{L}$.
The power flow on the line $ij$, $f_{ij}(\omega)$, is a linear function of the power injections $p(\omega) = g(\omega) - d(\omega)$:
\begin{equation} \label{eq:ptdf}
    f_{ij}(\omega) = \Phi_{(\cdot,ij)} p(\omega).
\end{equation}
The matrix $\Phi$ denotes the DC power transfer distribution factors (DC-PTDFs) \cite{wood2013}, with $\Phi_{(\cdot,ij)}$ referring to the row of $\Phi$ corresponding to the line $ij$. \dkm{We note that the \mbox{DC-PTDF} formulation of the DC power flow equations, i.e., the combination of~\eqref{eq:balance} and~\eqref{eq:ptdf}, implicitly ensures power balance at every bus in the system~\cite{BienCherHarn14}. The DC-PTDF formulation of the \mbox{DC OPF} problem used in this paper is equivalent to alternative formulations that explicitly include variables for the voltage angles at every bus.}

\subsubsection{Cost function} The generators have a quadratic cost function in terms of active power generation:
\begin{equation}
\cost(x) = \tfrac{1}{2}x^TMx + v^Tx + k_0,\label{eq: CostFunction}
\end{equation}
where $M$ is a diagonal matrix with non-negative entries. We minimize the expected generation cost $\E\left[\cost(g(\omega))\right]$. Substituting \eqref{eq: AffinePolicy} and taking the expectation in \eqref{eq: CostFunction}, we obtain
\begin{align*}
    &\E\left[\cost(g(\omega))\right] \\
    &\qquad = \E\left[ \tfrac{1}{2}(g - \beta \Omega)^TM(g - \beta \Omega) + v^T(g - \beta \Omega) + k_0 \right] \\
    &\qquad = \cost(g) + \E[\Omega]\left( g^T M \beta - v^T\beta \right) + \tfrac{1}{2}\E[\Omega^2]\beta^TM\beta \\
    &\qquad = \cost(g) + \tfrac{1}{2}\V(\Omega)\beta^TM\beta,
\end{align*}

\subsubsection{JCC-OPF} With these modelling considerations, we formulate the JCC-OPF as
\begin{subequations} \label{eq: CC_DC-OPF}
\begin{align}
    \label{eq: CC_DC-OPF_obj}
    &\min_{g,\beta} \quad \cost(g) + \tfrac{1}{2}\V(\Omega)\beta^T M \beta \\
    &\st \quad \sum_{i \in \G} g_i - \sum_{i \in \B} d_i = 0, \label{eq: NomPowBal} \\
    &\hspace{27 pt} \sum_{i \in \G} \beta_i = 1, \label{eq: FlucPowBal}\\
    & \proba\left( \begin{array}{rl}
        f_{ij}^{\LB} \leq \Phi p(\omega) \leq f_{ij}^{\UB},  &\forall ij \in \lag \\
        g_i^{\LB} \leq g_i - \beta_i \Omega \leq g_i^{\UB}, &\forall i \in \G
\end{array} \right) \geq 1-\alpha. \label{eq: ProbCons}
\end{align}
\end{subequations}
The objective~\eqref{eq: CC_DC-OPF_obj} minimizes the expected cost. The deterministic constraints \eqref{eq: NomPowBal}, \eqref{eq: FlucPowBal} ensure power balance. The JCC \eqref{eq: ProbCons} enforces bounds on the line flows and generator outputs, $f_{ij}^{\LB},f_{ij}^{\UB}$ and $g_{i}^{\LB},g_{i}^{\UB}$, with  probability $1-\alpha$. Here, $\alpha$ represents the acceptable violation probability.

\section{Representation of Chance Constraints} \label{sec: Quantile as CC}

Constraint \eqref{eq: ProbCons} results in a conventional nonlinear inequality $\vartheta(g,\beta) \geq 1-\alpha$. Problems with smooth nonlinear inequalities can be efficiently solved if one can compute the values of the inequalities as well as their gradients. In~\cite{PenaLuedWaec19}, the authors propose a smooth sample-based approximation of chance constraints to efficiently solve chance-constrained problems. This section summarizes the method from \cite{PenaLuedWaec19} and discusses modifications needed to address \eqref{eq: CC_DC-OPF}. We refer to one realization of $\omega$ as a ``scenario'' and a set of scenarios as a ``sample''.


We begin by defining the $(1-\alpha)$-quantile of a \blue{generic} random variable $Y$, denoted by $Q^{1-\alpha}(Y)$:
$$Q^{1 - \alpha}(Y) = \inf \{ y \in \mathbb{R} \mid \proba(Y \leq y) \geq 1 - \alpha \}.$$
From the above definition, we know that $\proba(Y \leq 0) \geq 1 - \alpha$ is equivalent to $Q^{1 - \alpha}(Y) \leq 0$, where $Y$ is a random variable taking values in $\mathbb{R}$. This definition can be extended to random variables $Y \in \mathbb{R}^m$, $m > 1$, if we let $\widehat{Y} = \max_{j=1,\ldots,m} \{ Y_j \}$ and consider $Q^{1-\alpha}(\widehat{Y})$ instead.

We denote the probabilistic constraint \eqref{eq: ProbCons} as $\proba (c(g,\beta;\omega) \leq 0) \geq 1 - \alpha$, \blue{where the random constraint vector $c(g,\beta;\omega)$ is defined by}
$$c(g,\beta;\omega) = \begin{pmatrix}
\Phi p(\omega) - f^{\UB} \\
f^{\LB} - \Phi p(\omega) \\
g - \Omega \beta - g^{\UB} \\
g^{\LB} - g + \Omega \beta
\end{pmatrix}.$$
Then, $c(g,\beta;\omega) \in \mathbb{R}^m$, where $m = {2|\lag| + 2|\G|}$. \blue{Since constraints $c_j(g,\beta;\omega) \leq 0$ for $j = 1,\ldots,m$ are equivalent to $C(g,\beta;\omega) = \displaystyle\max_{j=1,\ldots,m} \{c_j(g,\beta;\omega)\} \leq 0$, the single chance constraint $\proba(C(g,\beta;\omega) \leq 0) \geq 1-\alpha$ is equivalent to the joint chance constraint $\proba(c_j(g,\beta;\omega) \leq 0,\, j = 1,\ldots,m) \geq 1-\alpha$}. Reformulating \eqref{eq: CC_DC-OPF} yields
\begin{align}
    \min_{g,\beta} &\quad \cost(g) + \tfrac{1}{2}\V(\Omega)\beta^T M \beta \label{eq: Quant_DC-OPF} \\
    \st &\quad Q^{1-\alpha}\left( C(g,\beta;\omega) \right) \leq 0, \nonumber\\
    & \nonumber \quad \text{Eqns.~\eqref{eq: NomPowBal}, \eqref{eq: FlucPowBal}.}
\end{align}

For continuous random variables, the $(1-\alpha)$-quantile is obtained by inverting the cumulative density function (cdf) at the $(1-\alpha)$-level. Thus, whenever $C(g,\beta;\omega)$ defines a continuous random variable for any fixed value of $(g,\beta)$, an approximation of the quantile can be obtained from an approximation of the cdf. 

\blue{For the rest of Section~\ref{sec: Quantile as CC}, we assume that $C(g,\beta;\omega)$ defines a smooth function. This is clearly not the case because $C(g,\beta;\omega)$ is the maximum of the linear constraints given by the vector $c(g,\beta;\omega)$. However, we postpone the discussion of the non-smoothness of $C(g,\beta;\omega)$ to Section~\ref{sec: SolAlg} in order to introduce the approximation of the quantile via the cdf in a simplified manner.}

One way to approximate the cdf is to consider a sample $\{\omega_1,\ldots,\omega_N\}$ of the random variable $\omega$. The empirical probability that the random variable $C(g,\beta;\omega)$ takes a value less than or equal to $t$ (i.e., the empirical cdf evaluated at $t$) is
\begin{align}
    F^N(t;g,\beta) = \tfrac{1}{N}\sum_{i = 1}^N \mathbbm{1}(C(g,\beta;\omega_i) \leq t), \label{eq: EmpiricalCDF}
\end{align}
where $\mathbbm{1}$ is the indicator function, i.e., $\mathbbm{1}(A)$ takes the value of~$1$ if $A$ occurs or zero otherwise. Note that $F^N$ is non-smooth since the indicator function is not continuous at zero. For $t=0$, $F^N$ is equivalent to the SAA approximation used in MIP approaches (see \cite{LuedAhme08}).

\begin{figure}[htbp]
\centering
\begin{subfigure}{0.25\textwidth}
  \centering
  \includegraphics[width=0.99\linewidth]{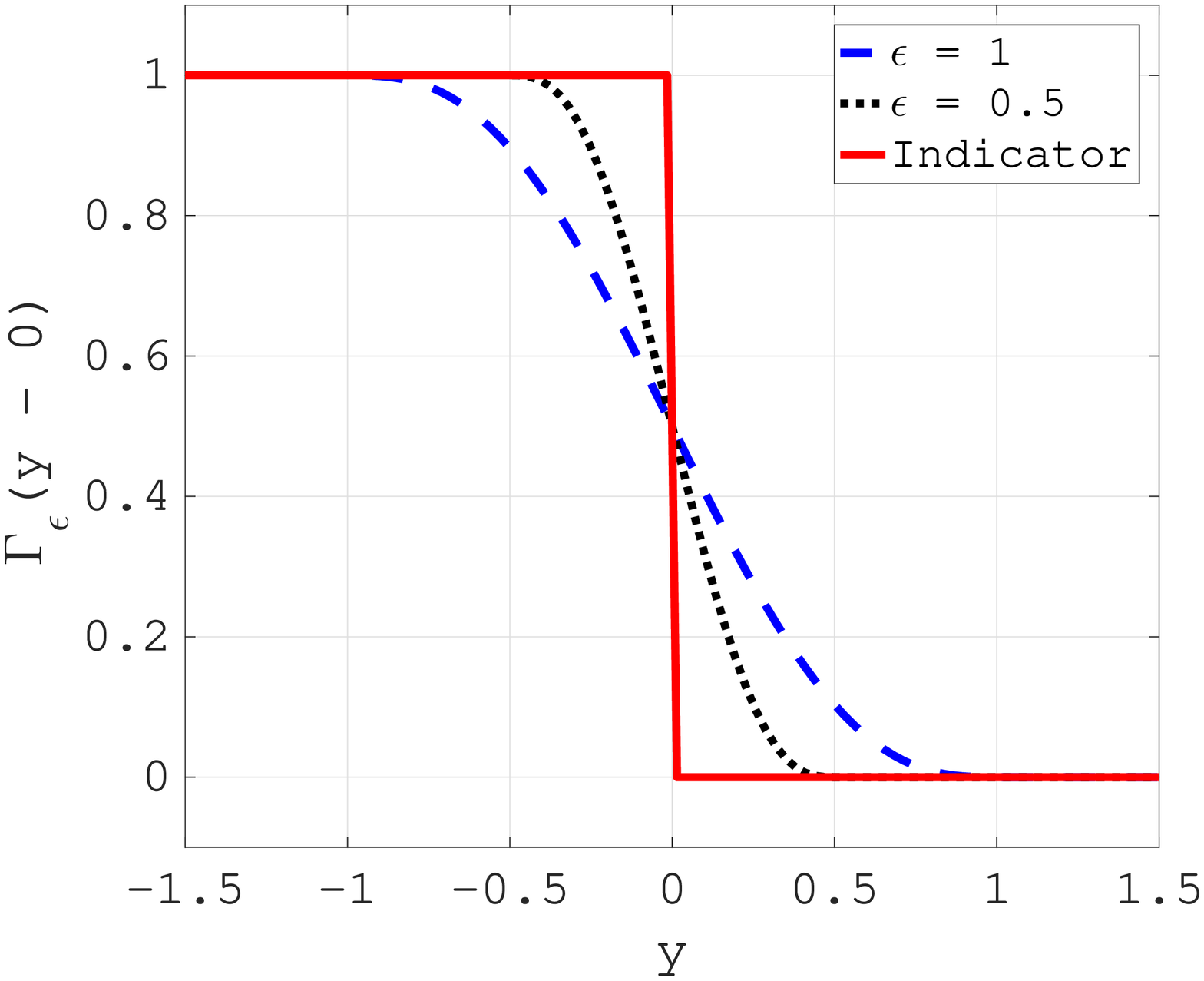}
  \caption{Different $\epsilon$ values.}
\end{subfigure}%
\begin{subfigure}{0.25\textwidth}
  \centering
  \includegraphics[width=0.99\linewidth]{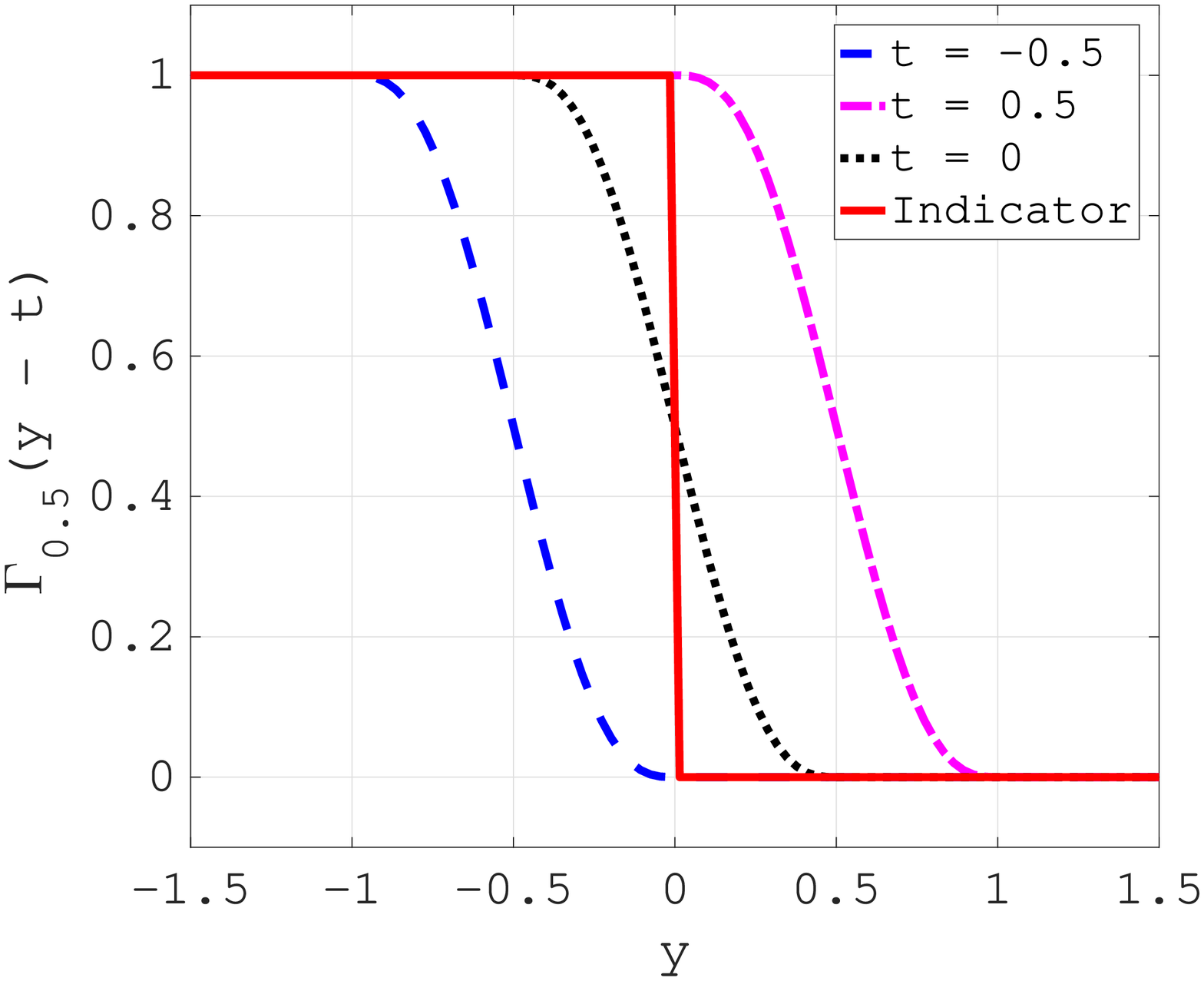}
  \caption{Different $t$ values.}
\end{subfigure}
\caption{Function $\Gamma_\epsilon(y-t)$.}
\label{fig: StepFunction}
\end{figure}

To obtain a smooth approximation of the cdf at the point $(g,\alpha)$, we follow an approach similar to \cite{GeleHoffKlop17,ShapDentRusz09} by defining
\begin{align}
    F_\epsilon^N(t;g,\beta) = \tfrac{1}{N}\sum_{i = 1}^N \Gamma_\epsilon (C(g,\beta;\omega_i)-t), \label{eq: CDF_SmoothAppr}
\end{align}
where $\epsilon > 0$ is a parameter of the following smooth approximation of the indicator function
\begin{align} \label{eq: SmoothIndicator}
\Gamma_\epsilon(y) = \left\{
	\begin{array}{lr}
		1,     & y \leq -\epsilon \\
  \gamma_\epsilon(y), & -\epsilon < y < \epsilon \\
		0,     &  y \geq \epsilon
	\end{array}
\right.
\end{align}
and $\gamma_\epsilon: [-\epsilon,\epsilon] \rightarrow [0,1]$ is a symmetric and strictly decreasing function such that $\Gamma_\epsilon$ is continuously differentiable. With this choice of $\gamma_\epsilon$, $F^N_\epsilon(t;g,\beta)$ is a differentiable approximation of the empirical cdf, $F^N(t;g,\beta)$ (see Fig.~\ref{fig: StepFunction}). We use the following $\gamma_\epsilon$ function based on the quartic kernel \cite[p. 353]{ScotTapiThom77}, which makes \eqref{eq: SmoothIndicator} twice continuously differentiable:
\begin{align}
    \gamma_\epsilon(y) = \frac{15}{16}\left( -\frac{1}{5}\left(\frac{y}{\epsilon}\right)^5 + \frac{2}{3}\left(\frac{y}{\epsilon}\right)^3 - \left(\frac{y}{\epsilon}\right) + \frac{8}{15} \right). \label{eq: gamma_eps}
\end{align}

For a fixed $g$ and $\beta$, the approximation of $Q^{1-\alpha}$ can be computed as the inverse of $F_\epsilon^N$ at $1-\alpha$. The inverse can be obtained from the value $Q_\epsilon$ such that
\begin{align}
    \sum_{i=1}^N\Gamma_\epsilon\left( C(g,\alpha;\omega_i) - Q_\epsilon \right) = N(1-\alpha). \label{eq: QuantRoot}
\end{align}
Reference~\cite{PenaLuedWaec19} shows that $Q_\epsilon$ results in an approximation of $Q^{1-\alpha}$ at $(g,\alpha)$. It also shows that the value $Q_\epsilon$ is unique, under mild conditions, and that it defines a function that maps the vector $C^N(g,\beta) = [C(g,\beta;\omega_1),\ldots,C(g,\beta;\omega_N)] \in \mathbb{R}^N$ to the root of \eqref{eq: QuantRoot}. We denote this function as $Q_\epsilon(C^N(g,\beta))$.

Hence, we propose the following approximation to \eqref{eq: CC_DC-OPF}:
\begin{subequations}\label{eq: ApproxProblem}
\begin{align}
    \min_{g,\beta} &\quad \cost(g) + \tfrac{1}{2}\V(\Omega)\beta^T M \beta \\
    \st &\quad \sum_{i \in \G} g_i - \sum_{i \in \B} d_i = 0 , \label{eq: Approx_1}\\
    &\quad \sum_{i \in \G} \beta_i = 1, \label{eq: Approx_2}\\
    &\quad Q_\epsilon \left( C^N(g,\beta) \right) \leq 0.\label{eq: QuantConstApprox}
\end{align}
\end{subequations}
\blue{Notice that $Q_\epsilon$ has taken the place of the chance constraint~\eqref{eq: ProbCons} following the fomulation in \eqref{eq: Quant_DC-OPF}.} Reference~\cite{PenaLuedWaec19} discusses the convergence and feasibility of this approximation with respect to the solutions of the original problem \eqref{eq: CC_DC-OPF} with increasing sample size. Section \ref{sec: SelectingParameters} details the choice of $\epsilon$ and~$t$.

Since $C(g,\beta;\omega)$ is the maximum of the entries of the vector $c(g,\beta;\omega)$, \eqref{eq: QuantConstApprox} results in a non-smooth function; this raises algorithmic challenges. Section~\ref{sec: SolAlg} discusses how to solve~\eqref{eq: ApproxProblem}.

\section{Solution Algorithm} \label{sec: SolAlg}

Adopting from the approach in~\cite{PenaLuedWaec19}, this section proposes an algorithm for solving~\eqref{eq: ApproxProblem}. To avoid having a non-smooth constraint, namely~\eqref{eq: QuantConstApprox}, we first reformulate~\eqref{eq: ApproxProblem} as an equivalent unconstrained optimization problem in which the constraints are added to the objective function via terms that penalize infeasible solutions. To address the challenges arising from the resulting non-smooth objective, we then propose an iterative algorithm that approximates the non-smooth unconstrained problem with a smooth constrained problem at each step of the algorithm. Finally, to improve tractability, we propose two modifications of this smooth constrained problem that make standard solvers compute the updates faster.


\subsection{$\ell_1$-penalty function}

Let $\pi > 0$ be a penalty parameter and $[x]^+ = \max\{0,x\}$. We propose an $\ell_1$-penalty function in order to solve \eqref{eq: ApproxProblem}:
\begin{align}
\phi_\pi(g,\beta) &= \cost(g) + \tfrac{1}{2}\V(\Omega)\beta^T M \beta + \pi\|V(g,\beta)\|_1 \label{eq: Penalization}
\end{align}
where $$V(g,\beta)=\left(\sum_{i \in \G} g_i - \sum_{i \in \B} d_i, \sum_{i \in \G} \beta_i - 1 , \left[ Q_\epsilon(C^N(g,\beta)) \right]^+ \right)$$ is the vector of constraint violations.
%
As shown in \cite{PenaLuedWaec19}, $\phi_\pi(g,\beta)$ is an exact penalty function, meaning that if $(g^*,\beta^*)$ is a local minimizer of $\phi_\pi$ for $\pi > 0$ and is feasible for problem \eqref{eq: ApproxProblem}, then $(g^*,\beta^*)$ solves \eqref{eq: ApproxProblem} \cite[p.~299]{HiriLema96}. Moreover, Theorem~2.1 in \cite{CurtOver12} shows that, under standard assumptions, there exists \blue{$\pi^* > 0$ such that the minimization of \eqref{eq: Penalization} yields a solution for \eqref{eq: ApproxProblem} for all $\pi \geq \pi^*$. This property makes the performance of exact penalty methods less dependent on the strategy for updating the penalty parameter than other penalty methods \cite[p. 507]{NoceWrig06}.}

\subsection{Minimizing the $\ell_1$-penalty function} \label{sec: Sl1QP}

To minimize \eqref{eq: Penalization}, we propose an S$\ell_1$QP-type trust-region algorithm that solves a sequence of quadratic programs (QP). At each iteration $k$, $\phi_\pi$ is approximated with a piecewise quadratic function that depends on $(g, \beta)$ at the current iteration. The trust region determines a region of the search space around the current iterate where the quadratic model provides a good approximation of the penalty function $\phi_\pi$. \blue{In each iteration of the algorithm, a trial step is computed as the minimizer of the model within the trust region.  If sufficient progress is made, the trial step is accepted. Otherwise, the trust region radius is reduced and a new trial step is computed.}

\blue{To be able to prove global convergence of the algorithm, the model must approximate the penalty function to first order. The following non-standard piecewise quadratic model was developed in \cite{PenaLuedWaec19} specifically for problems such as \eqref{eq: ApproxProblem} using a chain-rule-type approach,}
\begin{align}
& m(g,\beta,H;\delta) = \cost(g) + \tfrac{1}{2}\V(\Omega)\beta^T M \beta + \nabla \cost(g)^T \delta_g \nonumber \\
&\quad + \V(\Omega)\beta^T M \delta_{\beta} + \tfrac{1}{2}\delta^T H \delta + \pi \left( \left| \sum_{i \in \G} (g_i + \delta_{g_i}) - \sum_{i \in \B} d_i \right| \right. \nonumber \\
&\quad \left. + \left| \sum_{i \in \G} (\beta_i + \delta_{\beta_i}) - 1 \right| + \left[ \widetilde{Q}_{\epsilon,C}(g,\beta;\delta) \right]^+ \right), \label{eq: QuadModel}
\end{align}
where $\delta = [\delta_g; \delta_\beta] \in \mathbb{R}^{2|\G|}$ \blue{represents the trial step taken from the current point $(g,\beta)$}, $H \in \mathbb{R}^{2|\G|\times 2|\G|}$ is a symmetric matrix, and
\begin{align*}
\widetilde{Q}_{\epsilon,C}(g,\beta;\delta) &= \widetilde{Q}_\epsilon(C^N(g,\beta);\widetilde{C}^N(g,\beta;\delta)-C^N(g,\beta)), \\
\widetilde{Q}_\epsilon(z;p) &= Q_\epsilon(z) + \nabla Q_\epsilon(z)^T p, \\
[\widetilde{C}^N]_i(g,\beta;\delta) &= \max_j\{c_j(g,\beta;\omega_i) + \nabla c_j(g,\beta;\omega_i)^T \delta \}.
\end{align*}

Each iteration $k$ finds a descent step, $\delta^k$, for $\phi_\pi$ by minimizing the model $m(g^k,\beta^k,H^k;\delta)$ within a radius $\Delta^k$ for a given $H^k$. Minimizing $m(g^k,\beta^k,H^k;\delta)$ is challenging due to the non-smoothness introduced by the absolute values and the $\max$ operators that measure the infeasibility of the constraints. Hence, \eqref{eq: QuadModel} is rewritten as a smooth constrained QP by introducing slack variables $u$, $v$, and $w$:
\begin{subequations}\label{eq: MinModelOpt}%
\begin{align}%
&\min_{\delta,z,u,v,w} \; \cost(g^k) + \tfrac{1}{2}\V(\Omega)(\beta^k)^T M \beta^k + \nabla \cost(g^k)^T \delta_g \nonumber \\
&\qquad + \V(\Omega)(\beta^k)^T M \delta_{\beta} + \tfrac{1}{2}\delta^T H^k \delta + \pi \left[ (u + v)^T \pmb{1}_2 + w \right] \\
&\st \; \sum_{i \in \G} (g_i^k + \delta_{g_i}) - \sum_{i \in \B} d_i = u_1 - v_1 \label{eq: Mult1}\\
    &\qquad \sum_{i \in \G} (\beta_i^k + \delta_{\beta_i}) = 1 + u_2 - v_2 \label{eq: Mult2}\\
    &\qquad c(g^k,\beta^k;\omega_i) + \nabla c(g^k,\beta^k;\omega_i)^T \delta \leq z_i \pmb{1}_m, \; \forall i \in [N] \label{eq: MuMultipliers}\\[-0.5em]
    &\qquad \nabla Q_\epsilon(C^N(g^k,\beta^k))^T(z-C^N(g^k,\beta^k)) \nonumber\\
    &\hspace{100pt} + Q_\epsilon(C^N(g^k,\beta^k)) \leq w, \label{eq: LambdaMultipliers}\\
&\qquad t,u,w \geq 0, \quad \|\delta\|_\infty \leq \Delta^k, \label{eq: TrustRegionConstraint} 
\end{align}
\end{subequations}
where $\pmb{1}_n$ is the length-$n$ vector of ones and $[N] = \{1,\ldots,N\}$. The slack variables $u$, $v$ and $w$ ensure feasibility of the linearization of \eqref{eq: Approx_1}, \eqref{eq: Approx_2} and \eqref{eq: QuantConstApprox}, \blue{given by \eqref{eq: Mult1}, \eqref{eq: Mult2} and \eqref{eq: LambdaMultipliers}, respectively.} The $z$ variable in \cref{eq: MuMultipliers} represents the maximum of the linearization of $c$, i.e., $z_i = [\widetilde{C}^N]_i(g,\beta;\delta)$. Thus, \eqref{eq: MinModelOpt} is indeed equivalent to minimizing $m(g^k,\beta^k,H^k;\delta)$ with the addition of the trust-region constraint \eqref{eq: TrustRegionConstraint}.

A step $\delta^k$ obtained from solving \eqref{eq: MinModelOpt} is accepted if it results in sufficient decrease of $\phi_\pi$, i.e., we move in the $\delta^k$ direction only if the value $\phi_\pi(g^k+\delta_g,\beta^k+\delta_\beta)$ is sufficiently smaller than $\phi_\pi(g^k,\delta^k)$. If the step is accepted, we update $g^{k+1} = g^k +\delta_g$, $\beta^{k+1} = \beta^k + \delta_\beta$  and choose $\Delta^{k+1} \geq \Delta^k$; otherwise, the iterates are not accepted and we choose $\Delta^{k+1} < \Delta^k$.

For fast local convergence, $H^k$ is chosen as
\begin{align}
    &H^k = H_{\E} \label{eq: HessJCCP} \\
    &+ \lambda^k \nabla \overline{C}^N(g^k,\beta^k) \left[\nabla^2 Q_\epsilon({C}^N(g^k,\beta^k)) \right] \left[ \nabla \overline{C}^N(g^k,\beta^k) \right]^T, \nonumber
\end{align}
where $H_{\E}$ is the Hessian of the expected cost given by
$$H_{\E} = \begin{pmatrix}
    M & 0 \\
    0 & \V(\Omega)M
\end{pmatrix},$$
\blue{$\nabla \overline{C}^N(g^k,\beta^k)$ represents the transpose of the Jacobian of a smooth approximation of $C^N(g^k,\beta^k)$ (see \cite{PenaLuedWaec19}) obtained from}
\begin{align*} 
    & [\nabla \overline{C}^N(g^k,\beta^k)]_{\boldsymbol{\cdot}i} = \nabla c(g^k,\beta^k;\omega_i)\bar{\mu}^k_i, \;  \forall i \in [N],\\
    & [\bar{\mu}^k_i]_j = \frac{[\mu^k_i]_j}{ \lambda^k \left[ \nabla Q_\epsilon(C^N(g^k,\beta^k)) \right]_i}, \; \forall \, j \in [m], \; \forall \, i \in [N],
\end{align*}
and $\lambda^k$ and $\mu_i^k$ are the multipliers corresponding to \eqref{eq: LambdaMultipliers} and \eqref{eq: MuMultipliers}, respectively, from the previous iteration. If $\lambda^k = 0$ or $\left[ \nabla Q_\epsilon(C^N(g^k,\beta^k)) \right]_i = 0$, select one $j$ such that $c_j((g^k,\beta^k;\omega_i)) = C(g^k,\beta^k;\omega_i)$ and define $[\bar{\mu}^k_i]_j = 1$ and $[\bar{\mu}^k_i]_{\ell} = 0$ if $\ell \neq j$.

Lastly, as the stopping criterion of the algorithm, we focus on the infinity norm of
\begin{align}
  &\nabla \bar{\mathcal{L}}(g^k,\beta^k,\nu^k,\lambda^k) = \nabla \cost_{\E}(g^k,\beta^k) + \nu_1^k e_{\G} + \nu_2^k e_{\beta} \nonumber \\
  &\qquad + \lambda^k \sum_{i=1}^N \left[ \nabla Q_\epsilon(C^N(g^k,\beta^k)) \right]_i \nabla c(g^k,\beta^k;\omega_i)\bar{\mu}_i^k, \label{eq: GradJCCP}
\end{align}
where $\bar{\mathcal L}$ is an appropriately chosen Lagrangian function and
$ \nabla \cost_{\E}(\cdot)$ represents the gradient of the expected cost,
$$ \nabla \cost_{\E}(g,\beta) = \begin{pmatrix}
    Mg + v \\
    \V(\Omega)M\beta
    \end{pmatrix};$$
$\nu_1$ and $\nu_2$ are the multipliers associated with \eqref{eq: Mult1} and \eqref{eq: Mult2}, respectively; $\lambda$ and $\bar{\mu}$ are defined as before; and $e_\mathcal{I} \in \mathbb{R}^{2|\G|}$ is a vector such that $[e_\mathcal{I}]_i = 1$ if $i \in \mathcal{I}$ and $0$ otherwise. The function $\bar{\mathcal{L}}$ approximates the Lagrangian of a smooth optimization problem whose KKT points coincide with KKT points of \eqref{eq: ApproxProblem} (see (5.10) in \cite{PenaLuedWaec19}). Thus, if $\|\nabla \bar{\mathcal{L}}(g^*,\beta^*,\nu^*,\lambda^*)\|$ is less than a small convergence tolerance, then the point $(g^*,\beta^*,\nu^*,\lambda^*)$ is returned as a stationary point of \eqref{eq: Penalization}. If constraints \eqref{eq: Approx_1}--\eqref{eq: QuantConstApprox} are satisfied by $(g^*,\beta^*)$, we conclude that $(g^*,\beta^*)$ is a stationary point of \eqref{eq: ApproxProblem}.

\subsection{Improving the computation time} \label{sec: CompPerformance}

\dkm{Directly applying the algorithm described thus far has limited tractability since the number of constraints imposed in large DC~OPF problems with many scenarios leads to computationally challenging instances. We next propose two extensions that improve the algorithm's computational scalability. The first extension is a Hessian approximation that enables application of faster convex QP solvers. The second extension is a lazy constraint generation technique.}

\subsubsection{\dkm{Convex Hessian approximation}}
Since $Q_\epsilon$ is a non-convex function, the associated Hessian matrix $H^k$ in \eqref{eq: HessJCCP} is not necessarily positive semi-definite. Thus, \eqref{eq: MinModelOpt} might not be convex. In general, non-convex QPs are more challenging to solve than convex QPs. Our experiments show that the times for finding a global minmizer with the non-convex QP solver in \texttt{CPLEX} are generally very large and that they increase with the number of buses. To improve tractability, we replace $H^k$ by the positive definite approximation described next.

First, notice that if $\nabla^2 Q_\epsilon({C}^N(g^k,\beta^k))$ is positive semi-definite, then \eqref{eq: HessJCCP} is positive semi-definite. Thus, we replace $\nabla^2 Q_\epsilon({C}^N(g^k,\beta^k))$ by a positive semi-definite approximation $\widehat{Q}^k$ in which all negative eigenvalues are replaced by zero \cite[Section 3.4]{NoceWrig06}. Let  $A\Lambda A^T$ denote the spectral decomposition of $\nabla^2 Q_\epsilon({C}^N(g^k,\beta^k))$. We define $\widehat{Q}^k$ as
$$\widehat{Q}^k = A\left( \Lambda + \text{diag}(\tau_i) \right) A^T,$$
where
$$\tau_i = \left\{
	\begin{array}{lr}
		0,     &  \lambda_i \geq 0 \\
    \lambda_i, & \lambda_i < 0
	\end{array}
\right.$$
and $\lambda_i$ represents the $i$th eigenvalue of $\nabla^2 Q_\epsilon({C}^N(g^k,\beta^k))$. When we replace $\nabla^2 Q_\epsilon$ by $\widehat{Q}^k$, the semi-definite approximation of \eqref{eq: HessJCCP} is
\begin{align}
    \widehat{H}^k = H_{\E} + \lambda^k \nabla \overline{C}^N(g^k,\beta^k) \widehat{Q}^k \left[ \nabla \overline{C}^N(g^k,\beta^k) \right]^T. \label{eq: ApprHess}
\end{align}
This modification ensures that $\widehat{H}^k$ is positive semi-definite at every iteration. Substituting $\widehat{H}^k$ in \eqref{eq: MinModelOpt} makes the optimization problem convex and hence easier to solve. Furthermore, according to Theorem 5.3 in \cite{PenaLuedWaec19}, any choice of $H^k$ that is symmetric and bounded results in Algorithm \ref{alg: JCCPSolver} converging to a stationary point of $\phi_\pi$. Hence, while the number of iterates that the algorithm performs might increase, the approximation $H^k$ proposed in \eqref{eq: ApprHess} does not affect the convergence of the algorithm to a stationary point. The computation time of the extra iterates is offset by the time saved at each iterate by solving a convex QP instead \blue{of a nonconvex model}.

\subsubsection{\dkm{Lazy constraint generation}}
\blue{To further improve tractability, we utilize a lazy constraint generation technique.  Motivated by the observation that only a small fraction of the inequality constraints in the QP \eqref{eq: MinModelOpt} are active at the optimal solution~\cite{molzahn_roald-redundant2019}, we first solve a version of \eqref{eq: MinModelOpt} that only includes those inequalities that are either infeasible or within a certain threshold of becoming infeasible for the solution of the deterministic problem with $\omega = 0$. In an iterative fashion, we check which of the original inequalities are violated by the optimal solution of the reduced QP, add those to the QP, and solve the augmented problem until all of the original constraints are satisfied. For the next instance of \eqref{eq: MinModelOpt}, we start with the most recent set of inequalities.} \dkm{This procedure results in a reduction of up to $57\%$ of the total computation time for the test cases we considered.}


To summarize, Algorithm~\ref{alg: JCCPSolver} describes our proposed approach for solving \eqref{eq: ApproxProblem} for a given smoothing parameter $\epsilon > 0$. \blue{In Algorithm~\ref{alg: JCCPSolver} we represent the standard trust-region parameters by $\hat{\Delta} > 0$, the maximum trust-region radius; $\Delta_0 \in (0,\hat{\Delta})$, the initial trust-region radius; $\eta \in (0,1)$, the actual reduction ratio; and $\tau_1 \in (0,1)$ and $\tau_2 > 1$, the contraction and expansion coefficients of the trust-region radius. For the experiments in Section~\ref{sec: NumRes}, we set the values of these parameters to standard values of trust-region algorithms (see, e.g., \cite[Chapter 4]{NoceWrig06}).} The tolerances $\kappa_1 > 0$ and $\kappa_2 > 0$ represent the numerical accuracy for which we consider the problem optimal and feasible, respectively. 

\begin{algorithm}[htbp]
    \textbf{Inputs:} $\pi > 0$ (penalty parameter); $\hat{\Delta} > 0$, $\Delta_0 \in (0,\hat{\Delta})$, $\eta \in (0,1)$, $\tau_1 \in (0,1)$, and $\tau_2 > 1$ such that $1/\tau_2 \leq \tau_1$ (trust region parameters); $\kappa_1 > 0$ and $\kappa_2 > 0$ (optimality and feasibility tolerance); $(g^0, \beta^0, \nu^0, \lambda^0, \bar{\mu}^0)$ (initial point and multipliers); set $k\gets0$
	\begin{algorithmic}[1]
	    \STATE Let $\mathcal{J} = \{ (j,i) \mid c_j(g^0,\beta^0;\omega_i) > -\kappa_2 \}$.
		\WHILE{$\|\nabla\bar{\mathcal{L}}(g^k,\beta^k,\nu^k,\lambda^k) \|_{\infty} > \kappa_1$ \textbf{or} $\|V(g^k,\beta^k)\|_\infty > \kappa_1$}
		\STATE Set $C^N(g^k,\beta^k;\omega_i) = \max_{j = 1,\ldots,m}\{c_j(g^k,\beta^k;\omega_i) \}$ for all scenarios, compute $Q_\epsilon(C^N(g^k,\beta^k))$, $\nabla Q_\epsilon(C^N(g^k,\beta^k))$, and $\widehat{H}^k$.
		\STATE Obtain $(\delta,\nu,\mu,\lambda)$ by solving \eqref{eq: MinModelOpt} with the constraints of the type \eqref{eq: MuMultipliers} given in $\mathcal{J}$ (\textbf{if} $\delta^k = 0$ \textbf{stop}, stationary point of $\phi_\pi$ reached).
		\WHILE{There exists $(j,i)$ such that $c_j(g^k+\delta_g,\beta^k+\delta_\beta;\omega_i) > -\kappa_2$}
		\STATE $\mathcal{J} = \mathcal{J} \cup \{ (j,i) \mid c_j(g^k+\delta_g,\beta^k+\delta_\beta;\omega_i) > -\kappa_2 \}$
		\STATE Resolve \eqref{eq: MinModelOpt} with the constraints of the type \eqref{eq: MuMultipliers} given in $\mathcal{J}$.
		\ENDWHILE
		\STATE Compute the ratio $\rho^k =  \frac{\phi_\pi(g^k,\beta^k) - \phi_\pi(g^k+\delta_g,\beta^k+\delta_\beta)}{ m(g^k,\beta^k,\widehat{H}^k;0) -  m(g^k,\beta^k,\widehat{H}^k;\delta^k)}$
		\IF {$\rho^k < \eta$}
    		\STATE $\Delta^{k+1} = \tau_1 \min\{\Delta^k, \|d^k\|_\infty\}$
     		\STATE $g^{k+1} = g^k$; $\beta^{k+1} = \beta^k$
     		\STATE $\nu^{k+1} = \nu^k$; $\lambda^{k+1} = \lambda^k$; $\bar{\mu}^{k+1} = \bar{\mu}^k$
		\ELSE
		\STATE $g^{k+1} = g^k + \delta_g$; $\beta^{k+1} = \beta^k + \delta_{\beta}$
		\STATE Set $\bar{\mu}_i = \frac{\mu_i}{ \lambda \left[ \nabla Q_\epsilon(C^N(g^k,\beta^k)) \right]_i}$, for all $i = [N]$.
		\STATE $\nu^{k+1} = \nu$; $\lambda^{k+1} = \lambda$; $\bar{\mu}^{k+1} = \bar{\mu}$
		\IF {$\rho^k \geq \eta$ \textbf{and} $\|d^k \|_\infty = \Delta^k$}
        	\STATE $\Delta^{k+1} = \min\{\tau_2\Delta^k,\hat{\Delta}\}$
    	\ELSE
    		\STATE $\Delta^{k+1} = \Delta^k$
     	\ENDIF
    	\ENDIF
    	\STATE $k = k+1$
		\ENDWHILE
	\end{algorithmic}
	\textbf{Return:} $(g^k, \beta^k, \nu^k, \lambda^k, \bar{\mu}^k)$, optimal solution and multipliers.
	\caption{S$\ell_1$QP trust-region algorithm for CC DC-OPF}\label{alg: JCCPSolver}
\end{algorithm}

\section{Selecting the Smooth-Quantile Parameters} \label{sec: SelectingParameters}

The smooth approximation of the quantile $Q_\epsilon$ is motivated by a kernel estimation of the cdf \blue{resulting in \eqref{eq: gamma_eps}} \cite[p. 256]{AlemBolaGuil12}. The properties of the kernel approximation of the cdf can be extended to those of the quantile~\cite{Azza81}. These properties imply that large values of the smoothing parameter, $\epsilon$, reduce the variance among the estimators obtained from different samples, but can lead to biased estimators that are either consistently infeasible or consistently conservative. 

\begin{figure}[htbp]
\centering
\begin{subfigure}{0.25\textwidth}
  \centering
  \includegraphics[width=0.99\linewidth]{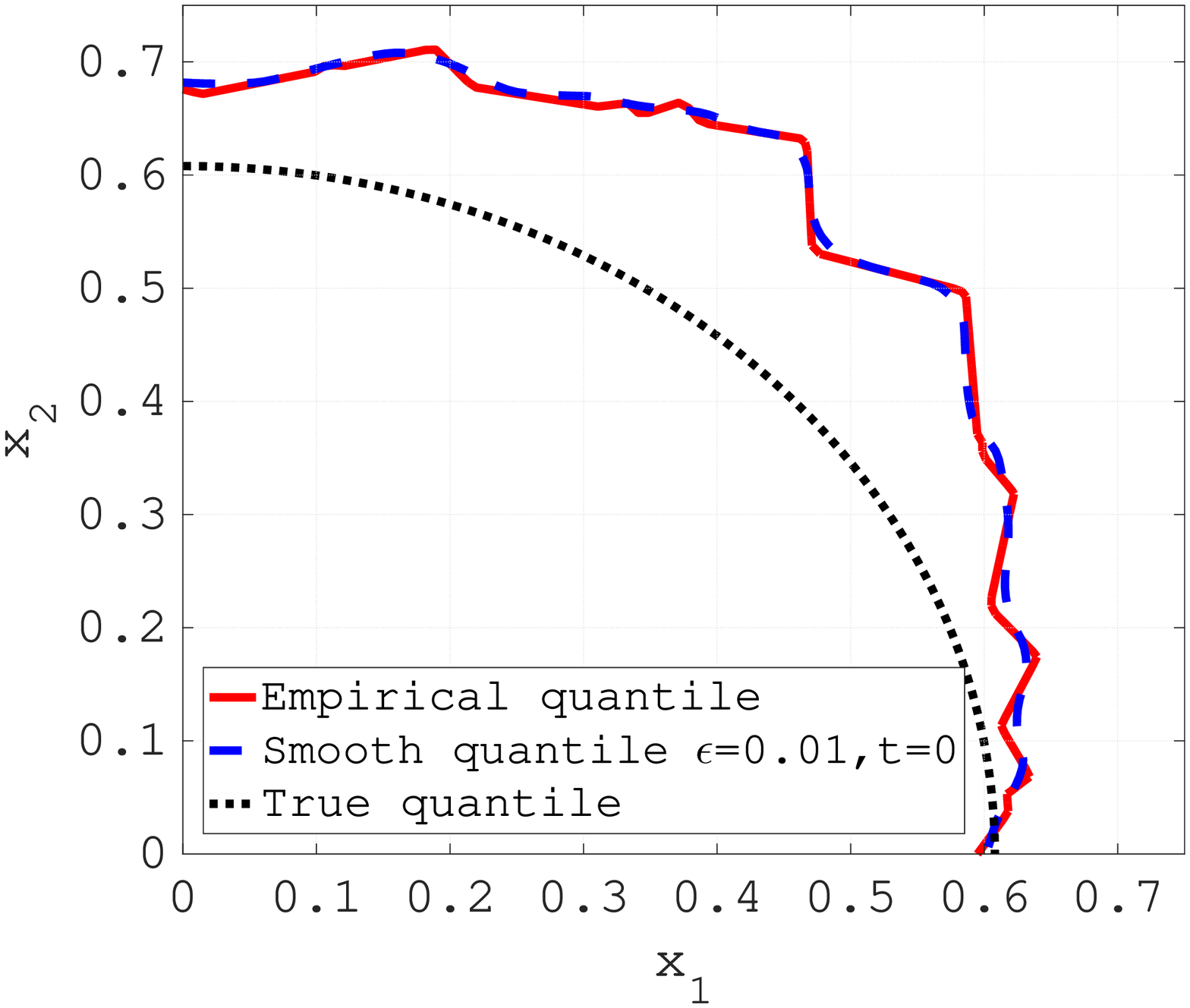}
\end{subfigure}%
\begin{subfigure}{0.25\textwidth}
  \centering
  \includegraphics[width=0.99\linewidth]{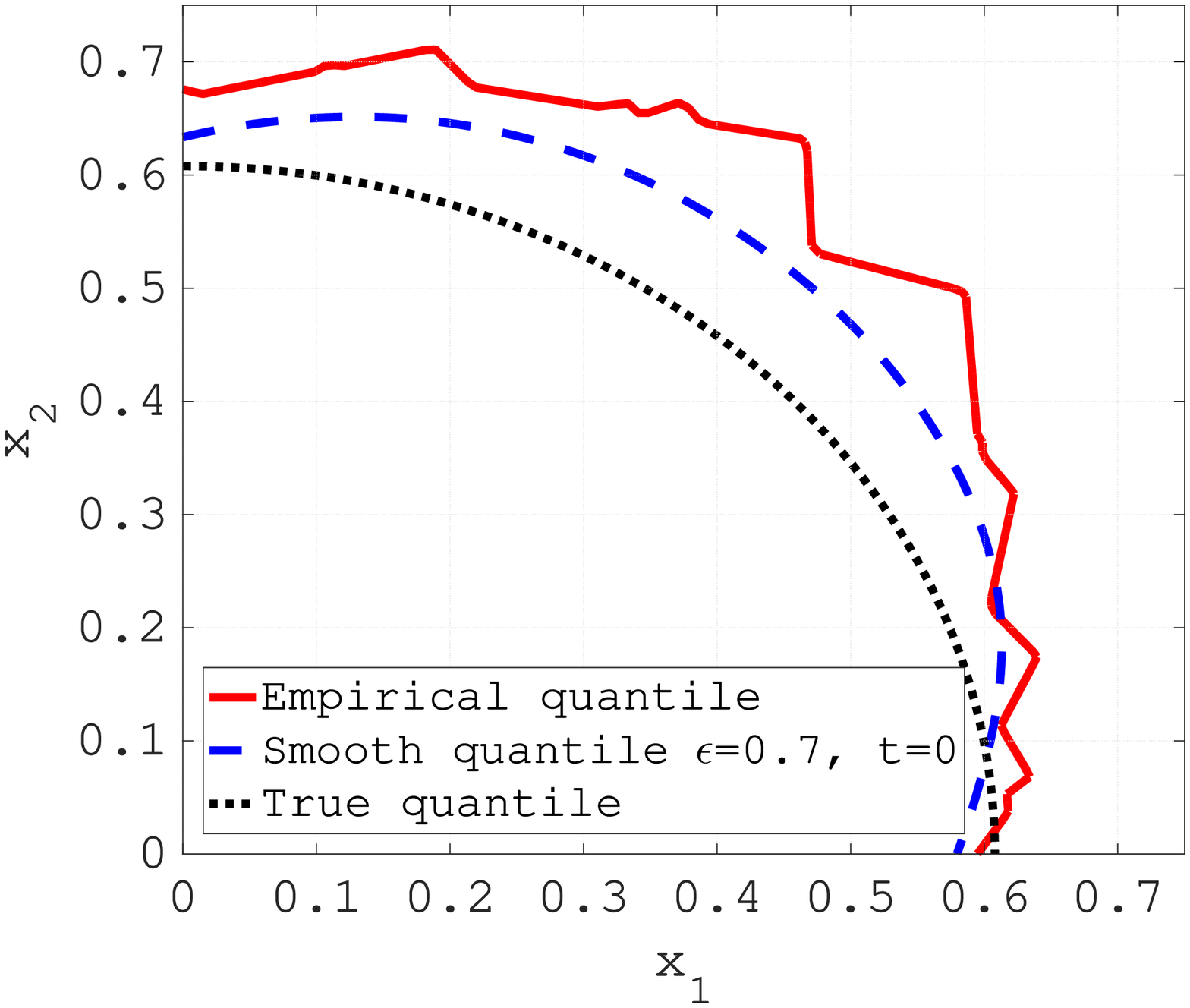}
\end{subfigure}
\begin{subfigure}{0.25\textwidth}
  \centering
  \includegraphics[width=0.99\linewidth]{NormalExample_h11-t0}
\end{subfigure}%
\begin{subfigure}{0.25\textwidth}
  \centering
  \includegraphics[width=0.99\linewidth]{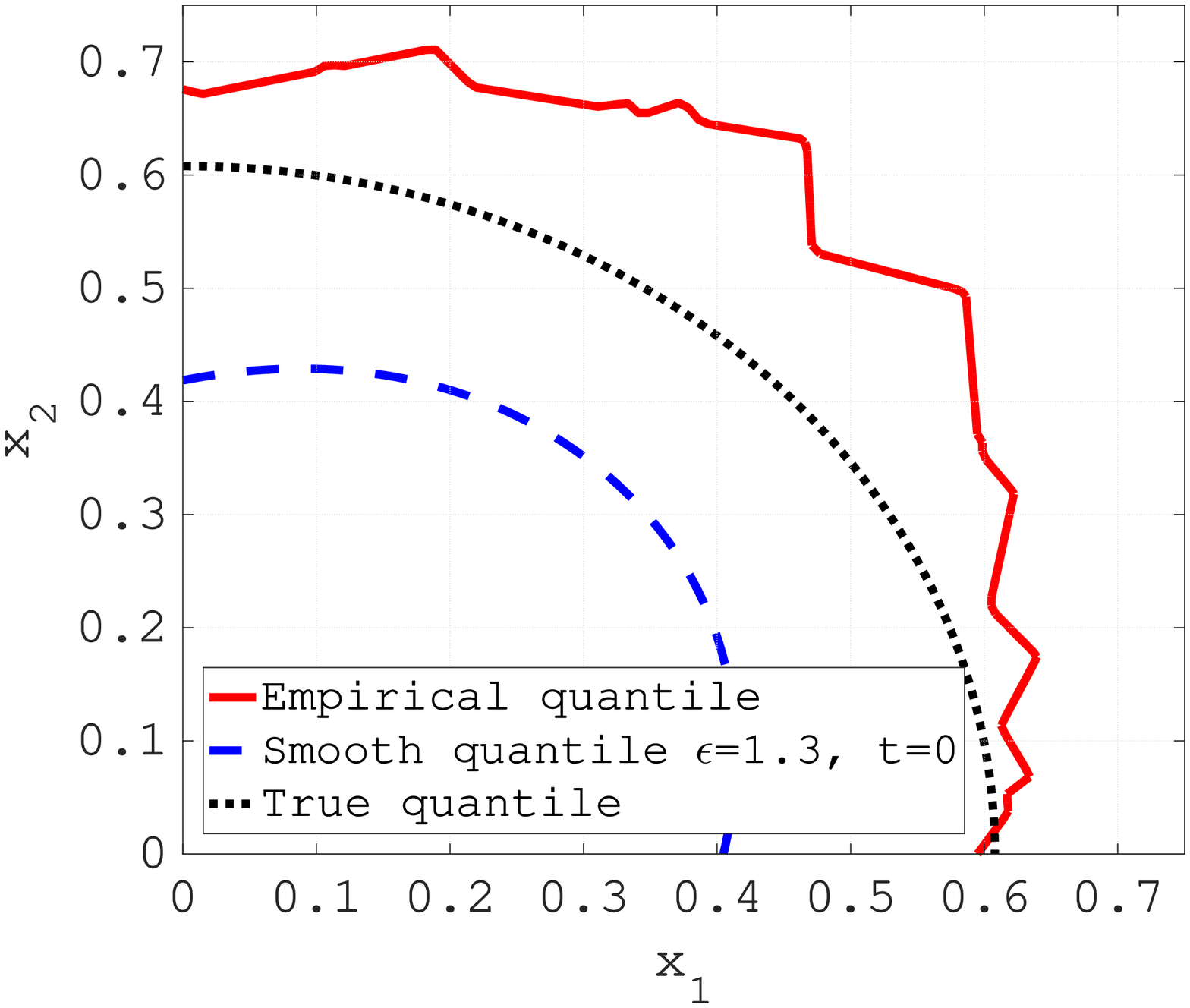}
\end{subfigure}
\caption{Comparison of different values of the smoothing parameter $\epsilon$.}
\label{fig: EpsilonTunning}
\end{figure}

\blue{To illustrate this, consider the following example in $\mathbb{R}^2$,
\begin{align}
    c(x_1,x_2;\omega_1,\omega_2) = x_1\omega_1 + x_2\omega_2 - 1, \label{eq: SimpleExample}
\end{align}
where $\omega_i \sim N(0,1)$ are independent random variables. The true feasible region of constraint $\proba\left(c(x_1,x_2;\omega_1,\omega_2) \leq 0 \right) \geq 0.95$ as well as its empirical and smooth approximations for a sample of size $N = 100$ are shown in Fig.~\ref{fig: EpsilonTunning}.}

\blue{Fig.~\ref{fig: EpsilonTunning} illustrates how the feasible region changes for different $\epsilon$. We observe that large values of $\epsilon$ introduce a bias in the feasible region of the smooth approximation by making it more conservative than the true feasible region. However,} increasing the value of $\epsilon$ also decreases the existence of spurious ``non-convexities'' in the feasible region. These ``non-convexities'' are not inherent to the problem but introduced by the discrete nature of the sample approximation, and they may cause local optimization algorithms (such as ours) to converge to local optima with worse objective values.
\blue{To avoid local optima and to increase the consistency of the solutions over different samples, using larger values of $\epsilon$ is advantageous. However, it is important to consider how the bias introduced by a large choice of $\epsilon$ can be counteracted.}

\begin{figure}[htbp]
\centering
\begin{subfigure}{0.25\textwidth}
  \centering
  \includegraphics[width=0.99\linewidth]{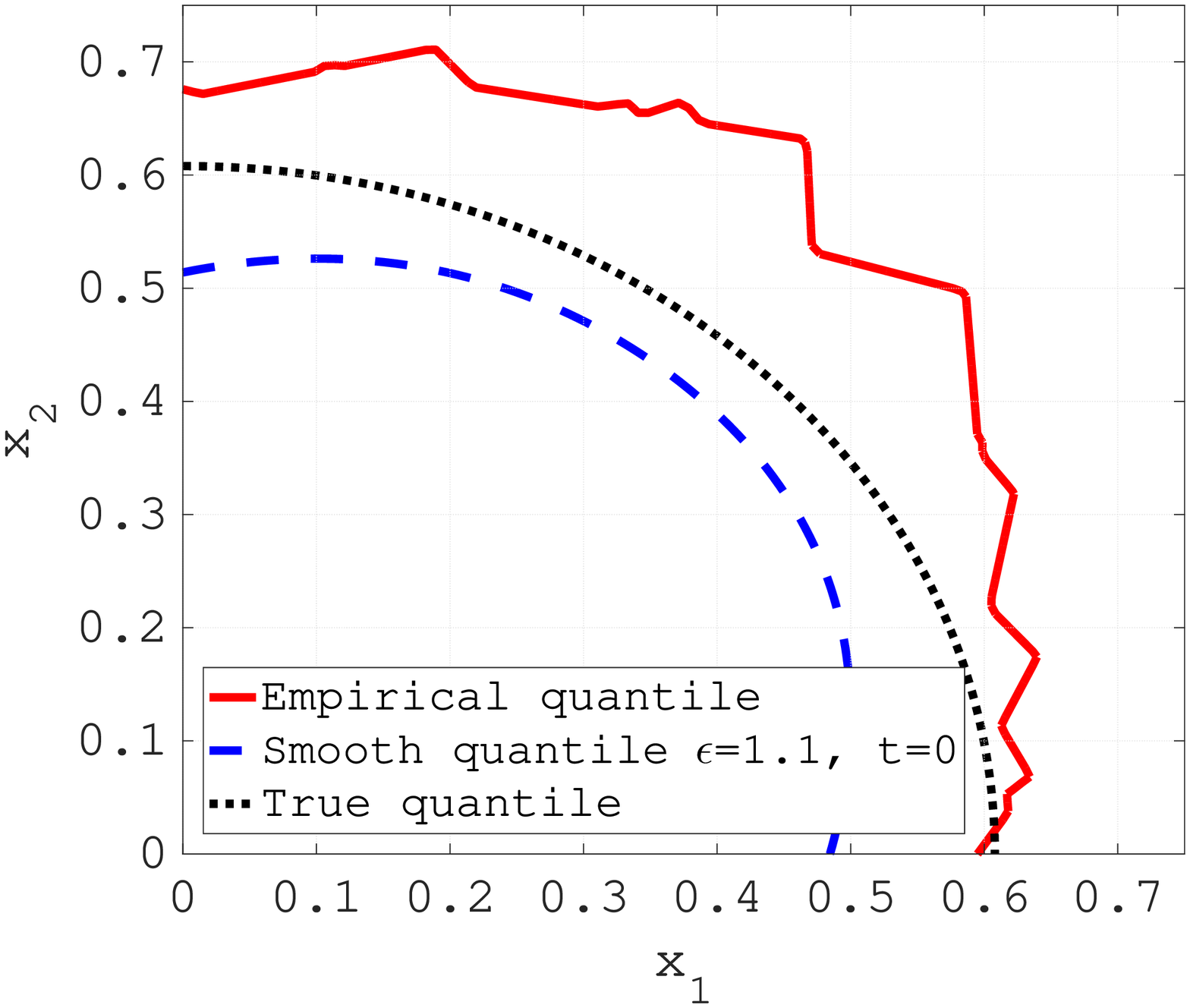}
\end{subfigure}%
\begin{subfigure}{0.25\textwidth}
  \centering
  \includegraphics[width=0.99\linewidth]{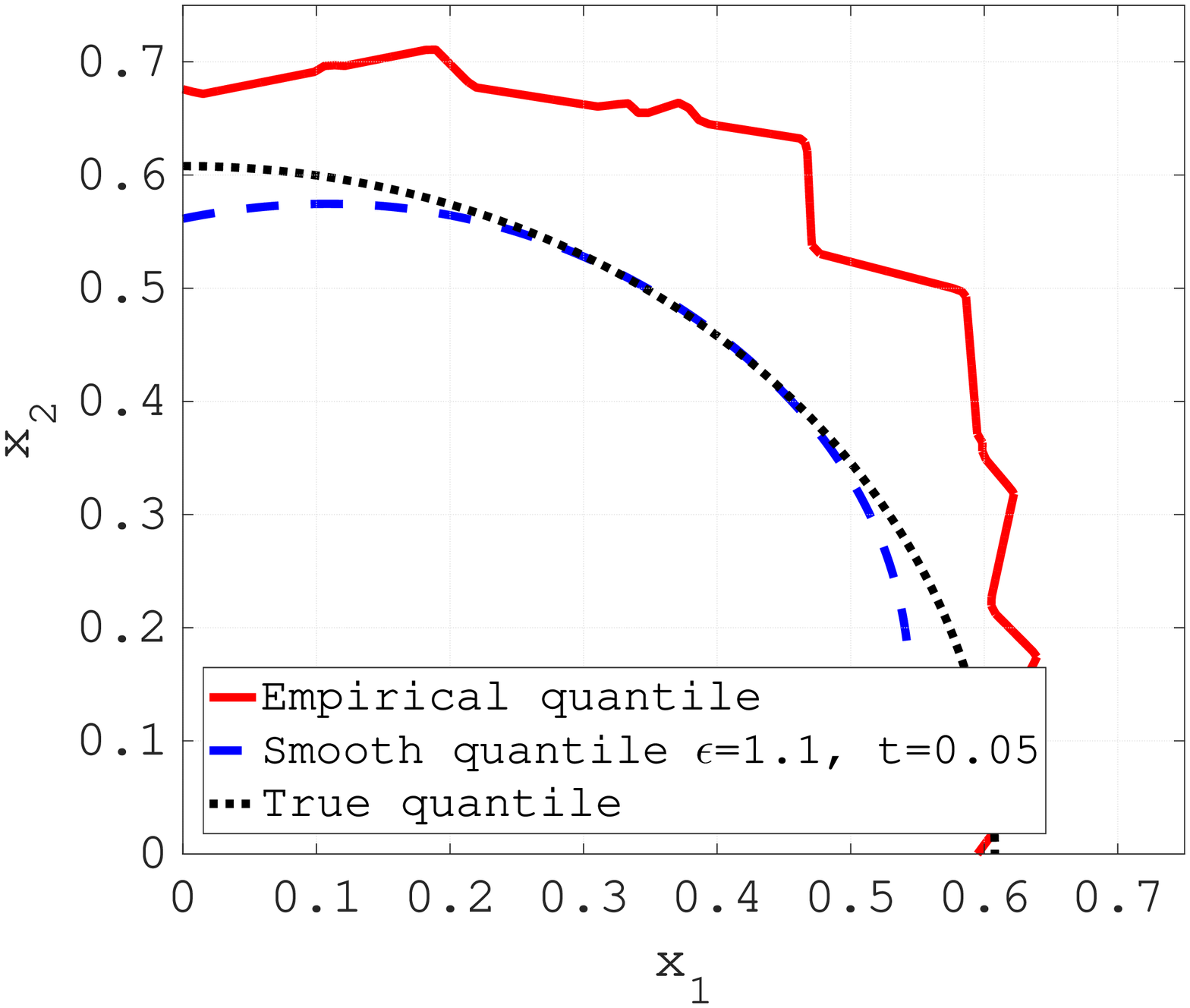}
\end{subfigure}
\caption{Comparison of different values of $t$.}
\label{fig: RHSTunning}
\end{figure}

To counteract the bias introduced by our choice of $\epsilon$, we propose to relax or strengthen the quantile constraint \eqref{eq: QuantConstApprox} by adjusting the right-hand side by $t \in \mathbb{R}$ as follows:
\begin{align}
    Q_\epsilon(C^N(g,\beta)) \leq t. \label{eq: QuantConst_Strong}
\end{align}
If $t < 0$, \eqref{eq: QuantConst_Strong} is more restrictive than \eqref{eq: QuantConstApprox}; if $t > 0$, \eqref{eq: QuantConstApprox} is relaxed. \blue{Hence, the feasible region of the approximated problem \cref{eq: ApproxProblem} gets smaller as the right-hand side~$t$ decreases (see Fig.~\ref{fig: RHSTunning}), which implies that the optimal objective value will increase as $t$ decreases. Given this monotone behavior of the approximated problem with respect to the right-hand side of \eqref{eq: QuantConst_Strong}, we propose a binary search method in order to find a value of $t$ such that the solution of \eqref{eq: ApproxProblem} attains an out-of-sample probability of $1-\alpha$ for a given $\epsilon$ and $N$.}

\blue{In Fig.~\ref{fig: RHSTunning}, it can be seen that for $N = 100$ and $\epsilon = 1.1$, the smooth quantile (blue dashed line) is slightly more conservative than the true quantile (black pointed line). However, if we keep $\epsilon = 1.1$ and adjust the right-hand side to $t = 0.05$, the smooth approximation is still conservative, in the sense that all feasible solutions for the approximated problem are feasible for the true problem, but now the bias is reduced.}

\blue{Our previous work in~\cite{PenaLuedWaec19} did not consider adjustments to~$t$. Here, we propose a new strategy that exploits the flexibility provided by the parameter~$t$ in combination with the parameter~$\epsilon$ in order to better tune the performance of our algorithm. We determine appropriate values for $\epsilon$ and $t$ in two stages. First, for a given problem class (i.e., a certain network), we empirically determine a suitable value of $\epsilon$ for a particular sample size. This value is then adjusted to different sample sizes using a statistical result obtained from the theory of kernel estimators \cite{Azza81} (see Section~\ref{sec: TuneEpsilon}). After this, for each particular instance, we tune the value of $t$ to obtain the best feasible point (see Section~\ref{sec: Tune t}}).



\subsection{Procedure for choosing the value of $\epsilon$} \label{sec: TuneEpsilon}

\blue{For a given sample size $\widehat{N}$, we want to select a value $\hat{\epsilon}$ that eliminates spurious ``non-convexities'' in the feasible region, while making sure that the bias introduced by this choice of the smoothing parameter is not too large (see Fig.~\ref{fig: EpsilonTunning}). We propose Algorithm~\ref{alg: EpsilonTunning} in order to select an appropriate $\hat{\epsilon}$ value.}

\begin{algorithm}[htbp]
    \textbf{Inputs:} $\widehat{N} \in \mathbb{N}$; $\hat{\epsilon}_0 > 0$; $\tau_{\epsilon} > 0$; $\widehat{M} \in \mathbb{N}$ \\ \vspace{-13pt}
	\begin{algorithmic}
	    \FOR{$m = 1,\ldots,\widehat{M}$}
	    \STATE Set $\epsilon_0 = \hat{\epsilon}_0$, $\epsilon_{LB} = 0$, $\epsilon_{UB} = \infty$, and $\ell = 0$.
		\REPEAT
		\STATE Obtain $(g^*,\beta^*)$ by solving \cref{eq: ApproxProblem} with $\widehat{N}$ samples, $t = 0$, and $\epsilon_\ell$.
		\STATE Approximate  $p_\ell \approx \proba(C(g^*,\beta^*;\omega) \leq 0)$.
		\IF{$p_\ell > 1 - \alpha$}
		\STATE $\epsilon_{\UB} = \epsilon_\ell$
		\STATE $\epsilon_{\ell+1} = (\epsilon_{\LB} + \epsilon_\ell)/2$
		\ELSE
		\STATE $\epsilon_{\LB} = \epsilon_\ell$
		\STATE \textbf{if} $\epsilon_{\UB} = \infty$ \textbf{then} $\epsilon_{\ell+1} = 2\epsilon_\ell$
		\STATE \textbf{else} $\epsilon_{\ell+1} = (\epsilon_{\UB} + \epsilon_\ell)/2$ \textbf{end if}
		\ENDIF
		\UNTIL{$|p_\ell - (1 - \alpha)| \leq \tau_{\epsilon}$ \textbf{or} $\epsilon_{\UB} - \epsilon_{\LB} \leq \tau_\epsilon$}
		\STATE $\bar{\epsilon}_m = \epsilon_\ell$
		\ENDFOR
		\STATE $\hat{\epsilon} = \max_{m=1,\ldots,\widehat{M}}\{\bar{\epsilon}_m\}$
	\end{algorithmic}
	\textbf{Return:} $\hat{\epsilon}$.
	\caption{Binary search for the smoothing parameter $\hat{\epsilon}$}\label{alg: EpsilonTunning}
\end{algorithm}

\blue{The goal of Algorithm~\ref{alg: EpsilonTunning} is to choose the most conservative approximation of the smoothing parameter $\epsilon$ from a set of $\widehat{M}$ samples. This value is obtained by running $\widehat{M}$ replications of the binary search algorithm described in \cite{PenaLuedWaec19} for a sample of size $\widehat{N}$ each. Then, since we have observed empirically that larger values of $\epsilon$ result in a more conservative approximation, we choose the maximum value of $\epsilon$ observed from the replications above in order to obtain the most conservative approximation of this parameter. In this paper, we choose $\widehat{N} = 100$ scenarios per sample and $\widehat{M} = 10$ samples.}

\blue{To avoid repeating Algorithm~\ref{alg: EpsilonTunning} for different sample sizes~$N$, we use a} result from \cite{Azza81} that approximates $\epsilon$ for different sample sizes. Reference \cite{Azza81} proves that, asymptotically, the optimal choice of $\epsilon$ that minimizes the mean square error between the true quantile and the approximated quantile is $\mathcal{O}(N^{-1/3})$. This can be used to estimate appropriate smoothing parameter values for sample sizes other than \blue{$\widehat{N}$. Thus, for a given sample size $N$, we consider $\epsilon = \tfrac{(\widehat{N}^{1/3})\hat{\epsilon}}{N^{1/3}}$.}

\subsection{Binary search to determine $t$} \label{sec: Tune t}

After choosing the smoothing parameter $\epsilon$ based on $N$, \blue{as described in Section \ref{sec: TuneEpsilon}}, we tune the parameter $t$ using the binary search described in Algorithm~\ref{alg: RHSTunning} for each separate sample. The parameter $\tau_t > 0$ in this algorithm determines the maximum difference allowed between the probability attained by the solution $(g^*,\beta^*)$ and the target probability $1-\alpha$.

\begin{algorithm}[htbp]
    \textbf{Inputs:} $t_0 = 0$; $t_{LB} = -\infty$; $t_{UB} = \infty$; $\tau_t,\kappa_t > 0$; $\ell = 0$ \\ \vspace{-15pt}
	\begin{algorithmic}
		\REPEAT
		\STATE Obtain $(g^*,\beta^*)$ by solving \cref{eq: ApproxProblem} with $t_\ell$.
		\STATE Approximate  $p_\ell \approx \proba(C(g^*,\beta^*;\omega) \leq 0)$.
		\IF{$p_\ell > 1 - \alpha$}
		\STATE $t_{\UB} = t_\ell$
		\STATE \textbf{if} $t_{\LB} = -\infty$ \textbf{then} $t_{\ell+1} = t_\ell -\kappa_t$ \textbf{end if}
		\STATE \textbf{if} $t_{\LB} > -\infty$ \textbf{then} $t_{\ell+1} = (t_{\LB} + t_\ell)/2$ \textbf{end if}
		\ELSE
		\STATE $t_{\LB} = t_\ell$
		\STATE \textbf{if} $t_{\UB} = \infty$ \textbf{then} $t_{\ell+1} = t_\ell + \kappa_t$ \textbf{end if}
		\STATE \textbf{if} $t_{\UB} < \infty$ \textbf{then} $t_{\ell+1} = (t_{\UB} + t_\ell)/2$ \textbf{end if}
		\ENDIF
		\UNTIL{$|p_\ell - (1 - \alpha)| \leq \tau_t$ \textbf{or} $t_{\UB} - t_{\LB} \leq \tau_t$}
	\end{algorithmic}
	\textbf{Return:} $t_\ell$, $(g^*,\beta^*)$, and $p_\ell$.
	\caption{Binary search for the right-hand side $t$}\label{alg: RHSTunning}
\end{algorithm}

When implementing Algorithm \ref{alg: RHSTunning}, we use the optimal solution and multipliers obtained when solving for the right-hand side $t^\ell$ as the initial points and multipliers for solving the problem with $t^{\ell+1}$. Once the binary search terminates, we select the solution from the $\ell$th right-hand side iteration that is feasible and has the best objective value. This ensures that we select the best of all the considered values of $t$.

\section{Case Study} \label{sec: NumRes}

This section demonstrates our method (denoted as the ``NLP approach'') using variants of the IEEE 14-, 57-, and 118-bus systems from \texttt{pglib-opf}~\cite{BabaEtAl19}. We compare our method against the scenario-approach \cite{CalaCamp05,CampGaraPran09,NemiShap06_SA} and the deterministic problem with $\omega = 0$; we refer to the solution of the deterministic problem as the ``nominal solution''. The obtained solutions are said to be good if they are: (1) consistent over different samples, (2) feasible for the true problem (evaluated with an out-of-sample test), and (3) low cost.

All computations were executed on \verb+Ubuntu 16.04+ with 256GB RAM and two \verb+Intel Xeon+ processors each with ten 3.10GHz cores. The algorithm is implemented in \verb+Matlab R2015b+, using \verb+CPLEX 12.6.3+ to solve the QP in \eqref{eq: MinModelOpt}. \blue{We set the \texttt{CPLEX} parameter \texttt{barrier.colnonzeros} to 1.} We use the parameters  $\pi = 10$, $\hat{\Delta} = 10^6$, $\Delta_0 = 1$, $\eta = 10^{-8}$, $\tau_1 = 1/2$, $\tau_2 = 2$, $\kappa_1 = 10^{-6}$, $\kappa_2 = 0.1$, $\tau_t=10^{-4}$ and $\kappa_t=0.01$.

To initiate the search for the right-hand side $t$, i.e., when $t^0 = 0$, we choose the initial points and multipliers for Algorithm \ref{alg: JCCPSolver} as follows: $g^0$ as the optimal solution of \eqref{eq: CC_DC-OPF} for $\omega = 0$, $\beta^0 = 1/|\G|$, $\lambda^0 = 0$, and $\bar{\mu}^0$ as described in Section~\ref{sec: SolAlg} using $\lambda^0 = 0$. For subsequent values of $t^\ell$, we initiate $g^0$, $\beta^0$, $\lambda^0$, and $\bar{\mu}^0$ using the optimal solutions returned from solving the problem with the previous right-hand side, $t^{\ell-1}$.

\vspace{-3pt}
\subsection{Uncertainty modeling}

The experiments in this section are based on normally distributed loads, i.e., $\omega \sim N(\Vec{0},\Sigma)$, where $\Sigma$ represents the covariance matrix.
\blue{To create the convariance matrix, we generated a $|\B| \times |\B|$ matrix, $A$, with entries taken from a uniform random variable with support in $[-1,1]$. Then, we obtained a positive definite matrix via $\widehat{A} = AA^T$. Finally, we scaled each entry of $\widehat{A}$ to ensure that $\Sigma_{ii} = \zeta d_i$ by defining
\vspace{-3pt}
$$\Sigma_{ij} = \zeta \tfrac{\widehat{A}_{ij}}{\sqrt{\widehat{A}_{ii}\widehat{A}_{jj}} } \left( \sqrt{d_i d_j} \right).$$
Here, $\zeta$ is a constant and $d$ is the vector of forecasted demands.} We let $\zeta = 0.1$ for cases 14 and 57; for case 118, we consider $\zeta = 0.05$ since the problem is infeasible for larger values of $\zeta$. In case 118, we also consider $\zeta = 0.01$ to compare the quality of the NLP solutions for different levels of variability. We aim to satisfy the probabilistic constraint at least $95\%$ of the time, i.e., $\alpha = 0.05$. For all solutions obtained in this section, out-of-sample approximations of $\proba(C(g^*,\beta^*;\omega) \leq 0)$ are computed using the empirical cdf with $N = 10^6$ scenarios.

\vspace{-3pt}
\subsection{Demonstration of joint chance constraints}

We first show the algorithm's performance with different sample sizes $N$. Table \ref{tab: OurMethod} presents the results of running 10 replications of the algorithm. The computation times are given in seconds; these times include the total time for the binary search algorithm to find the right-hand side $t$.

\begin{table}[htbp]
\begin{tabular}{c | cccc}
Case 14   & N \!=\! 100  & N \!=\! 200  & N \!=\! 500  & N \!=\! 1000 \\ \hline\hline
Min. obj (\$) &     2,106.2 &     2,120.1 &     2,104.6 &     2,105.4 \\ 
Avg. obj (\$) &     2,127.8 &     2,138.4 &     2,117.0 &     2,115.5 \\ 
Max. obj (\$) &     2,191.5 &     2,238.2 &     2,128.9 &     2,128.5 \\ \hline
Min. prob     &       0.950 &       0.950 &       0.950 &       0.950 \\ 
Avg. prob     &       0.950 &       0.950 &       0.950 &       0.950 \\ 
Max. prob     &       0.950 &       0.950 &       0.950 &       0.950 \\ \hline
Min. time (s) &      1.6812 &      4.6733 &      6.6583 &     28.4767 \\ 
Avg. time (s) &      3.2486 &      5.4729 &     15.2696 &     42.0556 \\ 
Max. time (s) &      6.5607 &      6.0644 &     35.0050 &     73.0396 \\ \hline
Avg. $t$ ($\times 10^{-3}$) & 17.219 & 8.3281 &  3.1875 &     -1.7344 \\
\multicolumn{1}{c}{} \\
Case 57       & N \!=\! 100  & N \!=\! 200  & N \!=\! 500  & N \!=\! 1000 \\ \hline\hline
Min. obj (\$) &      35,342 &      35,318 &      35,327 &      35,307 \\ 
Avg. obj (\$) &      35,413 &      35,394 &      35.358 &      35,337 \\ 
Max. obj (\$) &      35,508 &      35,471 &      35,381 &      35,357 \\ \hline
Min. prob     &       0.950 &       0.950 &       0.950 &       0.950 \\ 
Avg. prob     &       0.950 &       0.950 &       0.950 &       0.950 \\ 
Max. prob     &       0.950 &       0.950 &       0.950 &       0.950 \\ \hline
Min. time (s) &      2.5178 &      8.4458 &      11.586 &      64.078 \\ 
Avg. time (s) &      16.657 &      28.405 &      44.837 &      121.47 \\ 
Max. time (s) &      23.800 &      52.813 &     100.045 &      199.53 \\ \hline
Avg. $t$ ($\times 10^{-2}$) & 4.0375 & 2.7172 &  1.5625 &      0.9938 \\
\multicolumn{1}{c}{} \\
\!\!Case 118 ($\zeta \!=\! 0.01$) \!\! & N \!=\! 100  & N \!=\! 200  & N \!=\! 500  & N \!=\! 1000 \\ \hline\hline
Min. obj (\$) &     112,346 &     112,231 &     112,131 &     111,979 \\ 
Avg. obj (\$) &     112,594 &     112,469 &     112,231 &     112,038 \\ 
Max. obj (\$) &     112,891 &     112,626 &     112,453 &     112,160 \\ \hline
Min. prob     &       0.950 &       0.950 &       0.950 &       0.950 \\ 
Avg. prob     &       0.950 &       0.950 &       0.950 &       0.950 \\ 
Max. prob     &       0.950 &       0.950 &       0.950 &       0.950 \\ \hline
Min. time (s) &      13.143 &      45.256 &      268.08 &      658.34 \\ 
Avg. time (s) &      46.460 &      84.236 &      470.45 &      1127.1 \\ 
Max. time (s) &      86.017 &      181.70 &      879.82 &      1801.1 \\ \hline
Avg. $t$ ($\times 10^{-3}$) & 7.5234 & 3.9297 &  2.5156 &      4.7109 \\ 
\multicolumn{1}{c}{} \\
\!\!Case 118 ($\zeta \!=\! 0.05$)\!\!  & N \!=\! $100^*$ & N \!=\! 200 & N \!=\! 500 & N \!=\! 1000 \\ \hline\hline
Min. obj (\$) &     116,257 &     116,178 &     116,092 &     116,074 \\ 
Avg. obj (\$) &     116,615 &     116,315 &     116,138 &     116,107 \\ 
Max. obj (\$) &     117,183 &     116,670 &     116,168 &     116,165 \\ \hline
Min. prob     &       0.950 &       0.950 &       0.950 &       0.950 \\ 
Avg. prob     &       0.950 &       0.950 &       0.950 &       0.950 \\ 
Max. prob     &       0.950 &       0.950 &       0.950 &       0.950 \\ \hline
Min. time (s) &      26.842 &      104.34 &      581.31 &      1618.3 \\ 
Avg. time (s) &      50.336 &      142.75 &      867.83 &      2172.8 \\ 
Max. time (s) &      60.821 &      179.73 &      1925.9 &      3585.8 \\ \hline
Avg. $t$ ($\times 10^{-2}$) & 6.4642 & 5.1920 &  3.2102 &      2.0457 \\ 
\end{tabular}
\caption{Results from Algorithm \ref{alg: JCCPSolver} using $\hat{\epsilon}_{14} = 6.7(10)^{-2}$; $\hat{\epsilon}_{57} = 1.9(10)^{-1}$; $\hat{\epsilon}_{118} = 7.(10)^{-2}$ ($\zeta = 0.01$); $\hat{\epsilon}_{118} = 1.9(10)^{-1}$ ($\zeta = 0.05$). \textbf{*:} Statistics of feasible instances.} \label{tab: OurMethod}
\end{table}

Notice that the variability in the objective value decreases with increasing sample size. \blue{For example, in case 118 with $\zeta = 0.05$, the difference between the maximum and minimum costs decreases from 926 for sample size $N= 100$ to 76 for $N = 1,000$ scenarios, a $90\%$ decrease on the variability with respect to the sample size. However, even for $N = 100$, the variability between samples is no greater than $5.6\%$ for case 14 and $0.8\%$ for the rest of the cases.} This indicates that the NLP approach performs favourably even when using a small number of scenarios. Furthermore, the solutions returned for cases 14, 57 and 118 (with $\zeta = 0.01$) are always feasible for the true problem and \blue{accurately achieve the prescribed risk level of $95\%$.} For case 118 with $\zeta = 0.05$, all solutions are feasible for a sample size of at least 200. For $N = 100$, the solution obtained by the NLP approach on one instance is not feasible for the out-of-sample approximation of \eqref{eq: ProbCons}. \dkm{The achieved risk level is $92.1$\% for this instance.} We believe that this happens because the number of scenarios is insufficient for the level of variability.

\blue{For larger sample sizes, $N = 500$ and $N = 1,000$, the computation times can be improved by setting the \texttt{CPLEX} parameter $\texttt{barrier.colnonzeros} = 0$. This has shown a decrease of up to $55\%$ in the computation times. However, Tables~\ref{tab: OurMethod} and~\ref{tab: SA} only show the times obtained using $\texttt{barrier.colnonzeros} = 1$ for consistency in the paper.}

\vspace{-6pt}
\subsection{Comparison of the NLP and scenario approaches}

This section compares our solutions to those obtained from the scenario approach (SA) \cite{CalaCamp05}, which approximates \eqref{eq: CC_DC-OPF} as
\vspace{-2pt}
\begin{subequations} \label{eq: ScenAppr}
\begin{align}
    \min_{g,\alpha} &\quad c(g) \\
    \st & \quad f_{ij}^{\LB} \leq \Phi \hat{p}(\omega_s) \leq f_{ij}^{\UB}, \hspace{8 pt} \forall\, ij \in \lag, \; \forall s \in [N^{\text{SA}}]  \\
    &\quad g_i^{\LB} \leq g_i - \beta_i \Omega_s \leq g_i^{\UB}, \enskip \forall\, i \in \G, \; \forall s \in [N^{\text{SA}}]\\
    & \quad \text{Eqns.~\eqref{eq: NomPowBal}, \eqref{eq: FlucPowBal}.}
\end{align}
\end{subequations}
\vspace{-2pt}
where $N^{\text{SA}}$ is a pre-specified number of scenarios. SA specifies a minimum number of scenarios $N^{\text{SA}}$ such that a solution to \eqref{eq: ScenAppr}, which is feasible for all $N^{\text{SA}}$ scenarios, is also feasible for the probabilistic constraint \eqref{eq: ProbCons} with a probability of at least $1-\sigma$. \blue{An attractive feature of SA is that it results in a large-scale convex optimization problem, yet it tends to produce conservative results \cite{ZhanEtAl15}.}

We use the sample size given in \cite{CampGaraPran09} to select $N^{\text{SA}}$,
\vspace{-2pt}
$$N^{\text{SA}} \geq \tfrac{2}{\alpha}\left( \ln\left( \tfrac{1}{\sigma} \right) + n \right),$$
for $n = 2|\G|$, in this case. We select $\sigma = 10^{-4}$.
\vspace{-2pt}

Comparisons between the solutions obtained from the SA and NLP methods are presented in Table \ref{tab: SA}. The SA problem is solved for 10 different samples. We report the minimum, average, and maximum values of the objective and the out-of-sample probability of the returned solutions. If at least one of the instances is infeasible, the minimum probability is considered to be zero and the maximum objective function is marked as Inf. The average reported in the table does not consider the instances where the problem is infeasible. The number of infeasible instances using the SA approach for the different cases are: (1) Case 14: 7, (2) Case 57: 1, (3) Case 118 ($\zeta = 0.01$): 0, and (4) Case 118 ($\zeta = 0.05$): 3. \blue{We note that, since SA is a conservative approach, an infeasible SA instance does not indicate that the JCC problem is infeasible.}

First, observe that while the nominal solution is the least expensive, the solutions obtained when ignoring uncertainty are far from being feasible for cases 57 and 118; \blue{for these two cases, the nominal solution is feasible at most $38\%$ of the time.} In addition, there are many instances for which the SA algorithm cannot obtain a feasible solution to the approximated problem~\eqref{eq: ScenAppr}, \blue{even though the solutions obtained for the other samples clearly indicate that the problem is feasible. Of the 40 instances shown in the table, the SA approach was only able to return a solution for 29 of them. On the other hand, the NLP approach can always find a solution to the approximated problem \eqref{eq: ApproxProblem}, and out-of-sample testing verifies feasibility of the solutions to almost all of the problems we considered. The sole exception is one instance of case 118, with $\zeta = 0.05$ and $N=100$, that is infeasible with respect to the true problem; the out-of-sample risk level attained for this sample is $92.1\%$.}

Table~\ref{tab: SA} shows that the best solution from the SA can be up to $10\%$ more expensive than the worst NLP solution. \blue{For example, in case 14 the worst solution obtained by the NLP is $2,238.2$, while the best solution for the SA is $2,461.5$.} This demonstrates that the NLP approach provides solutions that are feasible without being overly conservative.

\blue{As several specific comparisons, the results for case 118 with $\zeta = 0.01$ and $N = 100$ show that the slowest computation time of the NLP method is faster than the fastest time of the SA, and that the average objective value from the NLP method is better. The results for case 118 with $\zeta = 0.05$ and $N = 100$ show that the best objective value from the SA is more expensive than the worst objective value from the NLP, and that the average computation time of the NLP is faster than that of the SA. Hence, there are instances for which the NLP method returns better solutions in less time than the SA.}

\blue{Finally, as the number of generators increases, the SA method prescribes a larger number of scenarios, $N^{SA}$, to be satisfied in order to guarantee feasibility. This significantly impacts the size of the SA problem and, as a consequence, the solution time of the SA method may become worse than the NLP (see case 118, $\zeta = 0.01$). For that reason, the advantages of the NLP method relative to the SA method are expected to be particularly pronounced for systems with many generators.}

\begin{table}[htbp]
\begin{tabular}{c | cccc}
Case 14       & Nominal  & SA        & NLP(100) & NLP(200) \\ \hline\hline
Min. obj (\$) &        - &   2,461.5 &  2,106.2 &     2,120.1 \\ 
Avg. obj (\$) &  2,051.5 &   2,489.0 &  2,127.8 &     2,138.4 \\ 
Max. obj (\$) &        - &      Inf  &  2,191.5 &     2,238.2 \\ \hline
Min. prob     &        - &         0 &    0.950 &       0.950 \\ 
Avg. prob     &    0.939 &     0.993 &    0.950 &       0.950 \\ 
Max. prob     &        - &     0.994 &    0.950 &       0.950 \\ \hline
Min. time (s) &        - &    0.0325 &   1.6812 &      4.6733 \\ 
Avg. time (s) &   0.0036 &    0.0425 &   3.2486 &      5.4729 \\ 
Max. time (s) &        - &    0.0492 &   6.5607 &      6.0644 \\ \hline
$N^{\text{SA}}$/Avg. $t$ & 1 &   516 &   0.0172 &      0.0083 \\
\multicolumn{1}{c}{} \\
Case 57        & Nominal    & SA          & NLP(100)     &  NLP(200)   \\ \hline\hline
Min. obj (\$) &           - &      35,493 &       35,342 &      35,318 \\ 
Avg. obj (\$) &      34,773 &      35,625 &       35,413 &      35,394 \\ 
Max. obj (\$) &           - &         Inf &       35,508 &      35,471 \\ \hline
Min. prob     &           - &           0 &        0.950 &       0.950 \\ 
Avg. prob     &       0.382 &       0.995 &        0.950 &       0.950 \\ 
Max. prob     &           - &       0.999 &        0.950 &       0.950 \\ \hline
Min. time (s) &           - &      0.2228 &       2.5178 &      8.4458 \\ 
Avg. time (s) &      0.0041 &      0.2776 &       16.657 &      28.405 \\ 
Max. time (s) &           - &      0.3761 &       23.800 &      52.813 \\ \hline
$N^{\text{SA}}$/Avg. $t$ & 1 &        637 &       0.0404 &      0.0272 \\
\multicolumn{1}{c}{} \\
Case 118 ($\zeta = 0.01$) & Nominal & SA & NLP(100) &    NLP(200) \\ \hline\hline
Min. obj (\$) &      - &         112,496 &  112,346 &     112,231 \\ 
Avg. obj (\$) & 109,791 &        112,824 &  112,594 &     112,469 \\ 
Max. obj (\$) &      - &         113,118 &  112,891 &     112,626 \\ \hline
Min. prob     &      - &           0.996 &    0.950 &       0.950 \\ 
Avg. prob     &  0.114 &           0.997 &    0.950 &       0.950 \\ 
Max. prob     &      - &           0.999 &    0.950 &       0.950 \\ \hline
Min. time (s) &      - &          86.022 &   13.143 &      45.256 \\ 
Avg. time (s) & 0.0071 &          88.546 &   46.460 &      84.236 \\ 
Max. time (s) &      - &          94.845 &   86.017 &      181.70 \\ \hline
$N^{\text{SA}}$/Avg. $t$ & 1 &      2998 &   0.0075 &      0.0039 \\
\multicolumn{1}{c}{} \\
Case 118 ($\zeta = 0.05$) & Nominal & SA & NLP(100)$^*$ & NLP(200) \\ \hline\hline
Min. obj (\$) &           - &    118,304 &      116,257 &     116,178 \\ 
Avg. obj (\$) &     109,791 &    122,553 &      116,615 &     116,315 \\ 
Max. obj (\$) &           - &        Inf &      117,183 &     116,670 \\ \hline
Min. prob     &           - &          0 &        0.950 &       0.950 \\ 
Avg. prob     &       0.057 &      0.996 &        0.950 &       0.950 \\ 
Max. prob     &           - &      0.997 &        0.950 &       0.950 \\ \hline
Min. time (s) &           - &     21.222 &       26.842 &      104.34 \\ 
Avg. time (s) &      0.0071 &     77.609 &       50.336 &      142.75 \\ 
Max. time (s) &           - &     94.166 &       60.821 &      179.73 \\ \hline
$N^{\text{SA}}$/Avg. $t$ & 1 &      2998 &       0.0646 &      0.0519 \\ 
\end{tabular}
\caption{\textbf{Nominal:} Solution of \eqref{eq: CC_DC-OPF} for $\omega = 0$. \textbf{SA:} SA with $\sigma = 10^{-4}$. \textbf{NLP(100):} Algorithm \ref{alg: JCCPSolver} with $N = 100$. \textbf{NLP(200):} Algorithm \ref{alg: JCCPSolver} with $N = 200$. \textbf{*:} Statistics of feasible instances.} \label{tab: SA}
\end{table}

\vspace{-6pt}
\section{Conclusions and Outlook} \label{sec: Conclusions}
This paper has developed a sample-based NLP algorithm for solving DC-OPF problems with JCC. By tuning two parameters in this algorithm using a proposed heuristic approach, the solutions obtained via this algorithm balance feasibility of the chance constraints and operational costs. Empirical results on several IEEE test cases demonstrate the algorithm's ability to jointly enforce chance constraints while being significantly less conservative with respect to operational costs than the alternative ``scenario approach'' proposed in prior literature. Our ongoing work is extending this approach to AC-OPF problems with JCC \dkm{as well as contingency constraints to model the possibility of component failures.}
\vspace{-6pt}

\section*{Acknowledgments}

The work by Alejandra Pe\~na-Ordieres was supported by the U.S. Department
of Energy, Office of Electricity Delivery and Energy Reliability under
contract DE-AC-02-06CH11357 and the National Science Foundation grant DMS-1522747.

The work by Daniel Molzahn was supported by the U.S. Department
of Energy, Office of Electricity Delivery and Energy Reliability under
contract DE-AC-02-06CH11357.

The work by Line Roald was supported by the Department of Energy, Office of Science, Office of Advanced Scientific Computing Research, Applied Mathematics program under contract DE-AC-02-06CH11347.

The work by Andreas W\"achter was supported by the National Science Foundation grant DMS-1522747.

\bibliographystyle{IEEEtran}
\bibliography{ref_ccacopf}

\end{document}